\documentclass{amsart}[12]

\makeatletter
\renewcommand{\@oddhead}{\thepage \hfil \textnormal{\small \scshape crepant transformation and the grassmannian flop} \hfil}
\makeatother

\usepackage{geometry}
\usepackage[dvipsnames]{xcolor}
\usepackage{amsmath}
\usepackage{amsthm}
\usepackage{mathrsfs}
\usepackage{amsfonts}
\usepackage{ amssymb }
\usepackage{amsmath,amscd}
\usepackage{array}
\usepackage{mathbbol}
\usepackage{amssymb}
\usepackage[all, cmtip]{xy}
\usepackage{tikz}
\usepackage[shortlabels]{enumitem}
\usepackage{tensor}
\usepackage{ stmaryrd }
\usetikzlibrary{arrows}
\usetikzlibrary{cd}
\usepackage{hyperref}
\usepackage{MnSymbol}

\renewcommand{\sslash}{\mathord{/\mkern-6mu/}}

\newtheorem{theorem}{Theorem}[subsection]
\newtheorem{lemma}[theorem]{Lemma}

\newtheorem{proposition}[theorem]{Proposition}
\newtheorem{definition}[theorem]{Definition}
\newtheorem{corollary}[theorem]{Corollary}
\newtheorem{remark}[theorem]{Remark}

\newtheorem{example}{Example}

\newcommand{\CC}{\mathbb{C}}
\newcommand{\QQ}{\mathbb{Q}}
\newcommand{\RR}{\mathbb{R}}
\newcommand{\PP}{\mathbb{P}}

\newcommand{\NN}{\mathbb{N}}
\newcommand{\ZZ}{\mathbb{Z}}
\renewcommand{\AA}{\mathbb{A}}
\newcommand{\UU}{\mathbb{U}}
\newcommand{\sgn}{\mathrm{sgn}}
\newcommand{\Alg}[1]{#1\mathrm{-Alg}}

\newcommand{\II}{\mathbb{I}}

\newcommand{\sO}{\mathscr{O}}

\newcommand{\bx}{\mathbf{x}}
\newcommand{\by}{\mathbf{y}}
\newcommand{\bt}{\mathbf{t}}
\newcommand{\bp}{\mathbf{p}}
\newcommand{\br}{\mathbf{r}}
\newcommand{\ba}{\mathbf{a}}

\newcommand{\cL}{\mathcal{L}}

\newcommand{\cY}{\mathcal{Y}}
\newcommand{\cH}{\mathcal{H}}
\newcommand{\cX}{\mathcal{X}}

\newcommand{\sL}{\mathscr{L}}

\newcommand{\Gm}{\mathbb{G}_m}
\newcommand{\bq}{\mathbf{q}}
\newcommand{\balpha}{\pmb{\alpha}}
\newcommand{\bbeta}{\pmb{\beta}}

\newcommand{\bH}{\mathbf{H}}
\newcommand{\rk}{\mathrm{rk}}

\newcommand{\bdr}{\mathbf{\dot{r}}}
\newcommand{\bdt}{\mathbf{\dot{t}}}
\newcommand{\bdp}{\mathbf{\dot{p}}}

\newcommand{\bs}{\mathbf{s}}

\newcommand{\cB}{\mathcal{B}}

\renewcommand{\sslash}{\mathord{/\mkern-6mu/}}
\newcommand{\one}{\mathbf{1}}
\newcommand{\bd}{\mathbf{d}}
\newcommand{\cM}{\mathcal{M}}
\newcommand{\cO}{\mathcal{O}}

\newcommand{\sH}{\mathscr{H}}

\newcommand{\bI}{\mathbf{I}}

\DeclareMathOperator{\Hom}{Hom}

\DeclareMathOperator{\Pic}{Pic}

\DeclareMathOperator{\Spec}{Spec}
\DeclareMathOperator{\Frac}{Frac}

\newenvironment{Comment}[2]{\noindent\color{#1}{\texttt #2: }}{}

\begin{document}



\title{Crepant Transformation Conjecture For the Grassmannian Flop}

\author{Wendelin Lutz}
\address{Department of Mathematics and Statistics\\
University of Massachusetts, Amherst\\
MA 01003\\
U.S.A}
\email{wendelinlutz@umass.edu}

\author{Qaasim Shafi}
\address{School of Mathematics\\
	Watson Building \\
	University of Birmingham\\
        Edgbaston
	B15 2TT\\
	UK}
\email{m.q.shafi@bham.ac.uk}

\author{Rachel Webb}
\address[R. Webb]{Department of Mathematics\\
Cornell University\\
Ithaca, NY 14853\\
U.S.A.}
\email{r.webb@cornell.edu}

\begin{abstract}We prove a highly explicit form of the Crepant Transformation Conjecture for Grassmannian flops. Our approach uses abelianization to first relate the restrictions of the Lagrangian cones to degree-2 classes, and then deduces the general result using ``explicit reconstruction'' (also known as the method of big $I$-functions). 
\end{abstract}

\maketitle
\setcounter{tocdepth}{1}
\tableofcontents

\section{Introduction}

% \Rachel{
% TODOS:
% \begin{itemize}
% %\item mention hiroshi's example 4.14 \url{https://arxiv.org/pdf/math/0411111.pdf}?
% \item any other acknowledgements?
% \item Mark's title: WALL CROSSING AND THE FOURIER–MUKAI TRANSFORM FOR
% GRASSMANN FLOPS
% \end{itemize}
% }

\subsection{Overview}
The problem of relating the Gromov-Witten invariants of $K$-equivalent birational varieties dates back to Ruan's crepant resolution conjecture \cite{ruan2006cohomology}. Here, we study the  Grassmannian flop introduced in \cite{DS14}, an example of two birational varieties $X\sslash_{\theta_+} G$ and $X\sslash_{\theta_-} G$ that are GIT quotients of a vector space $X$ by a general linear group $G$.
Very informally, our main theorem (Theorem \ref{thm:main}) states
\begin{equation}\label{eq:informal} \tag{\(\dagger\)}
\parbox{\dimexpr\linewidth-4em}{%
    \strut
    There are big points of the Lagrangian cones of $X\sslash_{\theta_\pm} G$ that are related by analytic continuation and a symplectic isomorphism of Givental spaces.%
    \strut
  }
\end{equation}
Before introducing the definitions necessary to make a precise statement, we make several remarks that should clarify the meaning and importance of the result.

First, a \textit{big point of the Lagrangian cone} of $X\sslash_{\theta_\pm} G$ is a series that determines the Lagrangian cone of $X\sslash_{\theta_\pm} G$ via a recursive algorithm known as Birkhoff factorization (see Definition \ref{def:big} and Proposition \ref{prop:bigI}). The Lagrangian cone in turn encodes all genus-zero descendant Gromov-Witten invariants of $X\sslash_{\theta_{\pm} }G$. In our case, the big points mentioned in \eqref{eq:informal} are explicit series in Chern classes on $X\sslash_{\theta_\pm} G$ and various formal variables.

Second, when $X$ is replaced by any representation of a torus $T$, for appropriate $\theta_+$ and $\theta_-$, the statement \eqref{eq:informal} is proved in \cite{coates2018crepant}. We realize our nonabelian wall crossing as a composition of the abelian wall crossings in \cite{coates2018crepant}; in particular, the analytic continuation in \eqref{eq:informal} is related via homotopy to a composition of analytic continuations in \cite{coates2018crepant}, and the symplectic isomorphism in \eqref{eq:informal} is the Weyl-anti-invariant part of a composition of the symplectic isomorphisms in \cite{coates2018crepant}.

Finally, an important idea in our proof is that while the objects most naturally related by analytic continuation are certain ``small $I$-functions'' of $X\sslash_{\theta_{\pm}} G$, we can recover an analytic continuation of big points via a technique called \textit{explicit reconstruction} (also known as the method of big $I$-functions; see \cite{giv15} and \cite{bigI}). While the use of reconstruction in crepant transformation problems is not new (it is used e.g. in \cite[Thm 3.12]{wc1}), as far as we know the observation that it can be applied to relate explicit points on the cone is new even in the toric case.

\subsection{Theorem statement}\label{sec:statement}
We recall the Grassmannian flop construction of \cite{DS14}.
Let $X = M_{k \times n} \times M_{n \times k}$ be the complex vector space consisting of pairs $(x, y)$ where $x$ is a $k \times n$ matrix and $y$ is a $n \times k$ matrix. Let $G = GL(k)$ act on this vector space by
\[
g\cdot(x,y) = (gx, yg^{-1}).
\]
%The integral characters of $G$ are powers of the determinant character, and the tensor product of this lattice with $\QQ$ is a 1-dimensional vector space. If we 
Define two characters $\theta_{\pm}$ of $G$ by 
\[\theta_+(g) = \det(g) \quad \quad \quad \quad \theta_-(g) = \det(g)^{-1}.\] 
The two GIT quotients $X\sslash_{\theta_\pm} G$ are both isomorphic to the direct sum of $n$ copies of the tautological bundle on $Gr(k, n)$, but their canonical morphisms to the stack quotient $[X/G]$ are different.

We work equivariantly with respect to the torus $K = (\CC^*)^{2n}$ whose action on $X$ is given by
\[
(k_1, k_2)\cdot(x,y) = (xk_1^{-1}, k_2y) \quad \quad \quad \quad (k_1, k_2) \in (\CC^*)^n \times (\CC^*)^n
\]
where we have identified an element of $(\CC^*)^n$ with a diagonal matrix. Note that this action commutes with the action of $G$ on $X$, hence descends to each of the quotients $X\sslash_{\theta_{\pm}}G$. 
Let $R_K$ denote the $K$-equivariant Chow ring of a point, tensored with $\QQ$, and define $F_K$ to be the field of fractions $\Frac(\CC \otimes R_K)$. Let $A^*_K(X\sslash_{\theta_{\pm}} G)_\QQ$ be the $K$-equivariant Chow group of $X\sslash_{\theta_\pm} G$ (tensored with $\QQ$). Let $N+1$ denote the rank of $A^*_K(X\sslash_{\theta_{\pm}} G)_\QQ$ as a free $R_K$-module. In the theorem below we will employ a formal variable $\log(\by^{\pm})$ associated to the unique hyperplane class and formal variables $\bx_1, \ldots, \bx_N$ implicitly associated to bases of the orthogonal complement to the span of the hyperplane class in $A^*_K(X\sslash_{\theta_{\pm}}G)_\QQ$. The associated Givental space is the $F_K((z^{-1}))$ vector space
\[H^{\pm}_G := F_K((z^{-1})) \otimes_{R_K} A^*_K(X\sslash_{\theta_{\pm}} G)_{\QQ}\]
equipped with a grading and a symplectic form (see Section \ref{sec:giv-module}).

Let $\PP^1$ have homogeneous coordinates $[\by^+: \by^-]$ and define
\[
p_+ := [0:1] \quad \quad \quad \text{and}\quad \quad \quad p_- := [1:0].
\]
Let $\widetilde{\PP^1}$ denote the formal analytic space whose underlying analytic space is the analytification of this $\PP^1$, and whose sheaf of functions is $\cO^{an}_{\PP^1} \otimes_{\CC} \CC\llbracket \bx_1, \bx_2, \ldots, \bx_{N}\rrbracket$. Finally, let $\Lambda^{\pm} := F_K\llbracket \bq^{\pm 1} \rrbracket$ denote the Novikov ring of $X\sslash_{\theta_{\pm}} G$.

\begin{theorem}\label{thm:main}
 There is a moduli space $\cM \simeq \widetilde{\PP^1}$, a degree-preserving symplectic isomorphism $\UU: H_G^+ \to H_G^-$ of $F_K((z^{-1}))$-modules, and series %\Qaasim{below should e.g. $\Lambda_+$ say $\Lambda^+_{F_K}$?}\Rachel{I changed the def of $\Lambda^{\pm}$ above. What do you think?}\Qaasim{Yep, that works!}
\[\begin{gathered}
\II_G^+ \quad \text{a big}\; \Lambda^+\llbracket \bx_1,  \ldots, \bx_{N}, \log(\by^+)\rrbracket\text{-valued point of the Lagrangian cone of $X\sslash_{\theta_+}G$}\\
\II_G^- \quad \text{a big}\; \Lambda^-\llbracket \bx_1, \ldots, \bx_{N}, \log(\by^-)\rrbracket\text{-valued point of the Lagrangian cone of $X\sslash_{\theta_-}G$}\\
\end{gathered}\]
such that 
\begin{itemize}
\item Setting $\bq=(-1)^{k-1}$ in the series $\II_G^{+}$ (resp. $\II_G^-$) yields a multivalued analytic function with values in $H_G^{+}$ (resp. $H_G^-$) on a punctured neighborhood of $p_+$ (resp. $p_-$), and $\UU {\II_G^+}|_{\bq=(-1)^{k-1}}$ analytically continues to $\II_G^-|_{\bq=(-1)^{k-1}}$ along a path $\gamma$.
\item Setting $\bq=1$ in the series $\II_G^{+}$ or $\II_G^-$ also yields a multivalued analytic function, and $\UU {\II_G^+}|_{\bq=1}$ analytically continues to $\II_G^-|_{\bq=1}$ along a different path $\delta$.
\end{itemize}
% Moreover, there is a tuple $(V, \widetilde{I^+}, \UU)$ where
% \begin{itemize}
% \item $V \subset \cM$ is a dense open set containing a path whose closure joins $p_+$ to $p_-$
% \item $\widetilde{I^+}$ is a multi-valued analytic function on $V$ with values in $H_+$
% \item $\UU: H^+ \xrightarrow{\sim} H^-$ is a degree-preserving symplectomorphism of $F_K((z^{-1}))$-modules
%\end{itemize}
% such that $\UU\widetilde{I^+}|_{\by^+ = (\by^-)^{-1}}$ (resp. $\widetilde{I^-}$) restricts to a branch of $I^+|_{\bq=1}$ (resp. $I^-|_{\bq=1}$) on the domain where $I^+$ (resp. $I^-$) is convergent.
\end{theorem}
 
The series $\II_G^{\pm}$ in this theorem are defined in \eqref{eq:computation} (but note the equivalent formula in \eqref{eq:def-IG}), and the meaning of the analytic continuation in the theorem statement is clarified in Remark \ref{rmk:meaning}. The matrix $\UU$ is obtained via an abelianization procedure from an analogous matrix arising from a Fourier-Moukai transform for a product of projective spaces. The reason for the two different specializations of the $\bq$ variable in Theorem \ref{thm:main}, leading to two different paths of analytic continuation, is that there are two different formulas relating the abelian and nonabelian $I$-functions (Theorem \ref{thm:nonabelianI}). 

Many of the arguments used to prove Theorem \ref{thm:main} generalize immediately to arbitrary representations $X$ of groups $G$ isogeneous to products of $SL(n)$'s (as far as we know, the restriction on the group is necessary for finding big points of the Lagrangian cone). One key place that we use the specific geometry of the Grassmannian flop is in the explicit homotopy of the path of analytic continuation in \eqref{eq:homotopy}.

\begin{remark}
A statement parallel to Theorem \ref{thm:main} was proved in \cite[Thm 6.1]{coates2018crepant} for arbitrary hyperplane wall crossings in toric variation of GIT. We remark that our Theorem does not include an analog of \cite[Thm 6.1(2)]{coates2018crepant} because the degree-2 cohomology of the affine quotient $X\sslash_0 G$ vanishes. Indeed, $X\sslash_0 G$ is the affine cone over a projective variety (as it is given by homogeneous equations), hence contractible in the analytic topology.

\end{remark}

%\Rachel{remark: for toric varieties, the cohomology ring is generated as an algebra by $H^2$, so in some sense it is ``enough'' to prove a statement about a ``small'' $I$-function that only parametrizes degree 2 classes. The Grassmannian chow ring is no longer generated in degree 2. However, it is \textit{equivariantly} generated in degree 2, and this is used in CITE frobman] .. . give more details/explanation. }

\begin{remark}\label{rmk:nonequivariant}
One can take the nonequivariant limit in Theorem \ref{thm:main}. However, the nonequivariant limits of the series $\mathbb{I}^{\pm}_G$ are independent of the $\log(\by^{\pm})$ variables and manifestly equal; i.e., the analytic continuation is just analytic continuation of a constant function and under some obvious choices of basis the isomorphism $\mathbb{U}$ is represented by the identity matrix. We work equivariantly in this paper in order to get a nontrivial result: in this setting, the relationship of $\mathbb{I}^+_G$ and $\mathbb{I}^-_G$ is not at all clear, and in the obvious bases the matrix $\UU$ is not the identity.
\end{remark}

\subsection{Relationship to other work}
Ruan's crepant resolution conjecture has been extended and refined by many authors, including Bryan-Graber \cite{bryan2005crepant}, Coates-Corti-Iritani-Tseng \cite{coates2009computing,wc1}, and Iritani \cite{iritani09}. For crepant transformations of smooth varieties, special cases and certain extensions of this conjecture have been proved for ordinary flops \cite{lee2016invarianceI,lee2016invarianceII,lee2016invarianceIII} and (as mentioned above) for flips arising from variations of geometric invariant theory (GIT) quotients by a torus \cite{coates2018crepant}.

Moreover, the work of Gonzalez-Woodward \cite{gonzalez2023wallcrossing} studies this conjecture for general variations of GIT quotient by a nonabelian group. Their work almost overlaps with our paper. However, it is not clear how the relationships found in \cite{gonzalez2023wallcrossing} compare to the ones presented here. Moreover, the theorem in \cite{gonzalez2023wallcrossing} is for projective GIT quotients of projective varieties, while our example concerns noncompact quotients of an affine variety. Finally, the theorem in \cite{gonzalez2023wallcrossing} holds nonequivariantly, while we work equivariantly---and in fact, as explained in Remark \ref{rmk:nonequivariant}, if one passes to the nonequivariant limit the example we study in this paper becomes trivial. 

When this paper was near completion, we learned of the parallel work of Priddis-Shoemaker-Wen \cite{PSW} which proves Theorem \ref{thm:main} in the case of the specialization $\bq=1$ (and in a relative setting). Our papers both use abelianization techniques and a homotopy argument. However, their paper constructs the matrix $\UU$ in Theorem \ref{thm:main} directly from a nonabelian Fourier-Moukai transform---let us call this matrix $\UU^{PSW}$---whereas our paper constructs $\UU$ as the anti-invariant part of the matrix coming from an abelian Fourier-Moukai transform. It is clear from the definition (see Proposition \ref{prop:UG-exists}) that our path $\delta$ in Theorem \ref{thm:main} agrees with the path of analytic continuation used in \cite{PSW}. Moreover, the proof of Lemma \ref{lem:unique} shows uniqueness of the matrix relating any two big points after analytic continuation along a given path. Hence we have
\[
\UU = \UU^{PSW}.
\] 
%\Rachel{I would like to say what this means for the FM transform and abelianization.}

% In a completely different vein, the paper \cite{BCFMV} is also closely related to the present work. In \cite{BCFMV} the authors show that the structure sheaf is a Fourier-Moukai kernel for the Grassmannian flop. One should be able to deduce our Theorem \ref{thm:main} from this result by identifying our matrix $\UU$ with this Fourier-Moukai transform via the Gamma-integral structure on the quantum cohomology of $X\sslash_{\theta_\pm} G$. (The analogous identification is used in \cite{coates2018crepant} to prove that the matrix $\UU$ appearing there is symplectic.) It would be interesting to explore this line of argument.

% This conjecture was extended and refined by many authors, including Bryan-Graber \cite{bryan2005crepant}, Coates-Corti-Iritani-Tseng \cite{coates2009computing,coates2009wall}, and Iritani CITE an integral structure].
% Special cases and certain extensions of this conjecture have been proved for ordinary flops \cite{lee2016invarianceI,lee2016invarianceII,lee2016invarianceIII}, for many examples of Hilbert-Chow morphisms CITE and for flips arising from variations of geometric invariant theory (GIT) quotients by a torus. \Rachel{find some more things to cite here}

\subsection{Summary of the paper}
Sections \ref{sec:givental}, \ref{sec:chow}, and \ref{sec:abelianization} apply to general GIT quotients of vector spaces. Section \ref{sec:givental} is a review of Givental formalism. Section \ref{sec:chow} contains results about the Chow ring of such quotients that we could not find in the literature, at least not in the generality needed for this project. Finally, Section \ref{sec:abelianization} is a mild restatement of the abelianization theorem in \cite{nonab1}. We specialize to the Grassmannian flop example and prove Theorem \ref{thm:main} in Section \ref{sec:prooftop}. Appendix \ref{sec:bigappendix} is again general and explains the notion of big points of the Lagrangian cone. The results in the appendix are not used to prove Theorem \ref{thm:main}, but they clarify the meaning and significance of our statement that the big points are related by analytic continuation.

\subsection{Acknowledgements}
The authors thank Tom Coates for suggesting the problem, including the example of the Grassmannian flop, as well as for helpful discussions. They especially thank Hiroshi Iritani for numerous enlightening conversations and most of the arguments in the Appendix. They are grateful to an anonymous referee for comments that clarified the exposition of the paper. The second author is supported by UKRI Future Leaders Fellowship through
grant number MR/T01783X/1. The third author thanks Victoria Hoskins and Mark Shoemaker for helpful conversations. The third author was partially supported by an NSF Postdoctoral Research Fellowship, award number 200213. 

\section{Background: equivariant Givental formalism for GIT quotients}\label{sec:givental}
In this section establish our definitions concerning the Givental space associated to an equivariant GIT datum.

Our input data is a tuple $(X, G, \theta, K)$ as follows.
The affine space $X$ is a representation of the connected complex reductive group $G$ and $\theta$ is a character of $G$. Define the loci of stable and semistable points $X^s_\theta(G)$ and $X^{ss}_\theta(G)$ as in \cite{king} and let $X\sslash_\theta G $ be the stack quotient $[X^{ss}_\theta(G)/G]$. We assume that $X\sslash_\theta G$ is a nonempty variety; recall that it is projective over the affine quotient $\mathrm{Spec}(\Gamma(X, \sO_X)^G)$, but may not be compact in general. Finally, $K$ is a complex torus with an action on $X$ that commutes with the $G$-action, so that the action of $K$ descends to $X\sslash_\theta G$. We assume that the fixed locus of the descended action is compact. Let $R_K$ be the equivariant Chow ring of a point (tensored with $\QQ$) and set $F_K = \Frac(\CC \otimes_\QQ R_K)$.

\subsection{Poincar\'e pairings and Gromov-Witten invariants}
To establish notation, we recall the definition of the $K$-equivariant Poincar\'e pairing for the potentially noncompact variety $X\sslash_\theta G$.
For a class $\gamma \in A^*_K(X\sslash_\theta G)_\QQ$, one defines
\begin{equation}\label{eq:def-integration}
\int_{[X\sslash G]_K} \gamma := \pi_*\sum_F \frac{\iota^*_F(\gamma)}{e_K(N_F)} \in \Frac(R_K)
\end{equation}
where the sum runs over the components $F$ of the fixed locus of $K$, and for each component $F$ we denote by $\iota_F$ its inclusion into $X\sslash _\theta G$ and by $e_K(N_F)$ the equivariant top Chern class of its normal bundle. The map $\pi_*$ is the equivariant pushforward to an equivariant point, tensored with $\Frac(R_K)$. By \cite[Cor 1]{EG98localization} the integral $\int_{[X\sslash G]_K} \gamma$ is the usual pushforward of $\gamma$ to an equivariant point when $X\sslash G$ is compact. We now recall the \textit{equivariant Poincar\'e pairing} 
\begin{align}
(\cdot, \cdot)_{X\sslash_\theta G} : A^*_K(X\sslash G)_\QQ \times A^*_K(X\sslash G)_\QQ &\longrightarrow \Frac(R_K) \label{eq:def-pairing}\\
(\alpha, \beta) &\longmapsto \int_{[X\sslash G]_K} \alpha \cup \beta.\notag
\end{align}
Though we will not need the formal definition in this paper, recall that the \textit{equivariant Gromov-Witten invariants} of $X\sslash_\theta G$ can be defined in an analogous way as localization residues (see e.g. \cite{liu12}).

\subsection{Givental space as a vector space}\label{sec:giv-module}
%\Rachel{concerned again: do I-funcs really have finite principal part wehn you set $q=1$? OK, yes and no: they can have infinite powers of $z$, but the coefficient of any $y$-monomial is polynomial in $z$, or at least has finitely many pos powers.}\Wendelin{Is this comment still relevant?}
The Givental space without Novikov variables\footnote{The Givental space without Novikov variables is the space where Gromov-Witten series related by analytic continuation live---see \cite[Thm 6.1]{coates2018crepant}.} is
\[
H := F_K((z^{-1})) \otimes_{R_K} A^*_K(X\sslash G)_\QQ,
\]
where $F_K((z^{-1}))$ denotes the ring of formal Laurent series in $z^{-1}$ with coefficients in the field $F_K$ and finite principal part. %\Wendelin{No Novikov here?}\Qaasim{Also not sure about why there is no novikov here}\Rachel{added a footnote} 
We view this as a graded $F_K((z^{-1}))$-vector space, with grading induced by the grading on the Chow group and by setting $\deg(z)=1$.
% There is a decomposition $H = H_+ \oplus H_-$ where
% \[
%  H_+ := F_K[z] \otimes_{R_K} A^*_K(X\sslash_\theta G)_\QQ \quad \quad \quad \text{and} \quad \quad \quad H_- := z^{-1}F_K\llbracket z^{-1}\rrbracket \otimes_{R_K} A^*_K(X\sslash_\theta G)_\QQ.
% \]

The Poincar\'e pairing on $A^*_K(X\sslash_\theta G)_\QQ$ induces a ``symplectic'' bilinear form $\Omega$ on $H$ given by
\[
\Omega(f(z), g(z)) = Res_{z=0}(f(-z), g(z))_{X\sslash_\theta G} \in F_K
\]
where by $(f(-z), g(z))_{X\sslash_\theta G}$ we mean to extend the Poincar\'e pairing on $A^*_K(X\sslash_\theta G)_\QQ$ to $H$ using $F_K((z^{-1}))$-linearity, and by $Res_{z=0}$ we mean to take the coefficient of $z^{-1}$.

\subsection{Givental space as a formal scheme}  There is a more general construction of the Givental space as a formal scheme which we now quickly summarize. We will not use this viewpoint much in the main body of the paper, but it will figure prominently in the Appendix, where more background will be given.

Let $\Lambda_{F_K}$ be the Novikov ring of $X\sslash_\theta G$ defined using numerical classes (rather than homological classes) of stable maps.\footnote{If $f:C \to X\sslash_\theta G$ is a stable map, its numerical class is the element of $\Hom(\chi(G), \ZZ)$ that sends a character $\zeta$ of $G$ to $\deg_C(f^*L_\zeta)$, where $L_\zeta$ is the line bundle on $X\sslash_\theta G$ induced by $\zeta$.} That is, if $\mathrm{Eff}(X\sslash_\theta G) \subset \Hom(\chi(G), \ZZ)$ is the monoid generated by numerical classes of stable maps, then $\Lambda_{F_K}$ is the completion of the monoid algebra $F_K[\mathrm{Eff}(X\sslash_\theta G)]$ with respect to the valuation
\[
q^\beta \mapsto \beta(\theta).
\]
The Novikov ring $\Lambda_{F_K}$ is a topological ring. Following \cite[Appendix B]{coates2009computing}, we define Givental's symplectic space $\sH$ as a formal scheme over $\Lambda_{F_K}$. We recall the full definition of $\sH$ in Definition \ref{def:givental-space}; for the main body of this paper, the reader can get by with the following facts.

\begin{itemize}
\item If $R$ is a topological $\Lambda_{F_K}$-algebra, an $R$-valued point of $\sH$ is an element of the module $R \otimes_{F_K} H$ (satisfying certain properties).
\item There is a formal subscheme $\sL \subset \sH$ called the \textit{Lagrangian cone of $X\sslash_\theta G$} that determines the genus-zero Gromov-Witten invariants of $X\sslash_\theta G$.
\end{itemize}

\section{Equivariant Chow groups of GIT quotients}\label{sec:chow}

In this section we relate the equivariant Chow groups of abelian and nonabelian quotients. For equivariant cohomology of compact quotients these statements are known and proofs can mostly be found in the literature (see for example \cite[Sec 3.1]{frobman}). We were unable to find references in the context of equivariant Chow rings of noncompact targets, and so for clarity we present a thorough exposition.

Section \ref{sec:notation} establishes notation on equivariant Chow groups that will be used throughout the paper. Sections \ref{sec:lim} and \ref{sec:W} are background, and for convenience, we state them in slightly more generality than the context of Section \ref{sec:notation}. The proofs of the results on Chow groups that we need are in Section \ref{sec:nonab chow ring}.

%The key results are two morphisms relating $A^*_K(X\sslash G)_\QQ$ and $A^*_K(X\sslash T)_\QQ$, one a ring morphism (with kernel) and the other a module isomorphism [ref], as well as an analog of Martin's celebrated integration theorem CITE martin thm B] leading to a relationship of the Poincar\'e pairings on $A^*_K(X\sslash G)_\QQ$ and $A^*_K(X\sslash T)_\QQ$ [ref].
\subsection{Notation}\label{sec:notation}
Let $(X, G, \theta, K)$ be as in Section \ref{sec:givental}, except that we do not require the fixed locus of $K$ to be compact.
We fix $T \subset G$ a maximal torus and we require that $X\sslash_\theta T$ is a nonempty variety (in particular the stable and semistable points agree). Let $W$ be the Weyl group of $T$ in $G$. Finally write $X\sslash_{\theta, G} T$ for the scheme quotient $X^{s}_\theta(G)/T$. There is a diagram of schemes
\begin{equation}\label{eq:geometry}
\begin{tikzcd}
X\sslash_{\theta, G} T \arrow[d, "g"] \arrow[r, "j"] & X\sslash_\theta T \\
X\sslash_\theta G
\end{tikzcd}
\end{equation}
where $j$ is an open embedding and $g$ is a $G/T$ fiber bundle. 

Fix a weight basis of $X$ with respect to the $T$-action. This basis determines a maximal torus of $GL(X)$ (equal to diagonal matrices), and we choose $K$ to be a subtorus of this maximal diagonal torus whose action on $X$ commutes with the $G$-action. The torus $T$ is also naturally embedded in the maximal diagonal torus of $GL(X)$ and it is possible that $T \cap K$ is nontrivial. 
%We write $T+K$ for the subgroup of diagonal matrices jointly generated by $T$ and $K$. Let $n$ be the rank of $K$. 
Any character $\chi$ of $T \times K$ induces a $K$-linearized line bundle on $X\sslash_{\theta} T$, yielding a homomorphism
\[
\chi(T \times K) \to \Pic^K(X\sslash_\theta T) \quad \quad \quad \quad \xi \mapsto \cL_\xi.
\]

% Any character $\chi$ of $T$
% induces a line bundle $\cL_\chi$ on $X\sslash_\theta T$. This yields the linearization map 
% \[
% \chi(T) \rightarrow \Pic(X\sslash_\theta T) \quad \quad \quad \quad \xi \mapsto \cL_\xi
% \]
% Similarly, a character of $T \times K$ induces a $K$-linearized line bundle on $X\sslash_{\theta, G}T$, yielding a linearization map  $\chi(K \times T) \rightarrow \Pic^K(X\sslash_\theta T).$

% If $\chi$ is a character of $T$ (resp. a character of $T \times K$
% \footnote{The group $T \times K$ is not the same as the subgroup of $GL(X)$ generated by $T$ and $K$; in fact, one can expect it to have higher rank.}) we use $\cL_\chi$ to denote the line bundle on $X\sslash_\theta T$ or $X\sslash_{\theta, G} T$ induced by $\chi$ (resp. the $K$-linearized line bundle induced by $\chi$). 

%We use $\xi_1, \ldots, \xi_N$ to denote the characters of $T \times K$ such that $(t, k) \cdot x_i = \xi_i(t, k) x_i$. Observe that the $\xi_i$ also define characters of $T$ or of $K$ by restriction. 

If $\cX$ is a smooth variety or algebraic global quotient stack of dimension $n$, we use $A^*(\cX)_\QQ$ to denote the group $A_*(\cX) \otimes_\ZZ \QQ$ together with its intersection product, where $A^{i}(\cX)_\QQ := A_{n-i}(\cX)_\QQ$. Recall that $A^*(\cX)_\QQ$ is isomorphic to the operational Chow ring \cite[Prop~4]{EG98}.

\subsection{Nonequivariant limit}\label{sec:lim}
%\Wendelin{Removed reference to a specific fiber, fixed $\PP^n$ vs $(\PP^n)^k$} \Rachel{Is it clear that the pullbacks are independent of choice of basepoint? I want to work with a fixed nonequivariant limit homomorphism (even if the particular one doesn't matter).}

Let $K$ be a complex torus as above, and let $\cX$ be a smooth separated variety with $K$-action. (For example, we might have $\cX = X\sslash_\theta G$ or $\cX = X\sslash_\theta T$.)
Let $N$ be a positive integer, let $\rk(K)$ denote the rank of $K$, and define $U_N = (\AA^{N+1} \setminus \{0\})^{\rk (K)}$. Recall that 
%if $\cX$ is a smooth separated variety %\Qaasim{is this supposed to say variety?}\Rachel{as opposed to?} \Qaasim{it's not really important, it's just we've used curlyX and a moment ago it says if $\cX$ is a smooth variety or algebraic global quotient stack.. and now it just says variety} \Rachel{Yes I don't know that this current sentence is true if $\cX$ is a global quotient artin stack. It is true if $\cX$ is DM (no need to be a global quotient even) but we don't need that.} 
%with $K$-action, 
we have (essentially by definition)
\[
A^{\leq N}_K(\cX)_\QQ = A^{\leq N}(\cX \times_K U_N)_\QQ.\]
 The projection $\cX \times_K U_N \to U_N/K$ is an $\cX$-fiber bundle over $U_N/K \cong (\PP^N)^{\rk (K)}$ 
%and choosing the point $p=[1:0:\ldots:0]$ in $U_N/K$ induces via 
and pullback to a fiber %over $p$ 
induces a graded ring homomorphism 
\[
A^{\leq N}_K(\cX)_\QQ \to A^{\leq N}(\cX)_\QQ
\]
independent of the choice of fiber.
%%%%%%%% REASON: I believe that pullback is Gysin pullback, i.e. intersection with the fiber, and it doesn't matter which fiber we intersect with because all fibers are rationally equivalent, because all points of $(\PP^N)^{\rk(K)}$ are rationally equivalent.
These pullbacks are compatible with the natural restrictions $A^{\leq N'}_K(\cX)_\QQ \to A^{\leq N}_K(\cX)_\QQ$ for $N' \geq N$, and we get a graded ring homomorphism
\begin{equation}\label{eq:limit}
A^*_K(\cX)_\QQ \to A^*(\cX)_\QQ
\end{equation}
which we call the \textit{nonequivariant limit}.
As usual, the Leray-Hirsch theorem lets us find splittings of \eqref{eq:limit}, as we explain in the following lemma.

%\Rachel{I removed all the tensor-Q's from this lemma, but maybe I shouldn't have?}
\begin{lemma}\label{lem:lh}
Let $\cX$ be a smooth proper variety with an action of a torus $K \times \Gm$, and assume that the restriction of the action to $\Gm$ has finitely many fixed points.
Let $\alpha_1, \ldots, \alpha_n \in A^*_K(\cX)$ be classes whose nonequivariant limits \eqref{eq:limit} are a $\ZZ$-basis for $A^*(\cX)$. Then the $\alpha_i$ freely generate $A^*_K(\cX)$ over the $K$-equivariant Chow ring of a point.
%\Rachel{Make sure anything using this satisfies properness.}

\end{lemma}

\begin{proof}
Let $N \geq \mathrm{max}_i(\deg(\alpha_i))$ be an integer. Recall that
\[
A_K^{\leq N}(\cX) = A^{\leq N}(\cX \times_K U_N).
\]
The variety $\cX \times_K U_N$ is an $\cX$-fiber bundle over the projective space $U_N/K$, and the $\alpha_i$ are classes in $A^*(\cX \times_K U_N)$ whose restrictions are a basis for the Chow group of a fiber. Moreover the fibration $\cX \times_K U_N \to U_N/K$ is trivial over each affine open subset of the base $U_N/K \simeq (\PP^N)^{\rk(K)}$ given by nonvanishing of a section of $\cO(1) \boxtimes \cO(1) \boxtimes \ldots \boxtimes \cO(1)$. Hence (using our assumption on the existence of the $\Gm$-action) we may apply the Leray-Hirsch theroem of \cite[Appendix C]{ColFul89}: an element of $A^*(\cX \times_K U_N)$ can be uniquely written as a linear combination of the classes $\alpha_i$ with coefficients in $A^*(U_N/K)$. By degree considerations, the unique expressions for classes in $A^{\leq N}(\cX \times_K U_N)$ have coefficients in $A^{\leq N}(U_N/K)$, which is just $R_{K, \ZZ}^{\leq N}$ by definition (we have used $R_{K, \ZZ}$ to denote the integral equivariant Chow ring of a point).

\end{proof}

\subsection{Weyl action}\label{sec:W} We begin with a general discussion of Weyl actions on vector spaces. Let $k$ be a field containing $\QQ$.
Recall the sign function $\sgn: W \to \{\pm 1\}$ defined by $\sgn(w) = \det(w)$ when $w$ is viewed as an automorphism of $\chi(T)$. A $W$-vector space is a $k$-vector space $H$ with a linear action by $W$. An element $h\in H$ is $W$-invariant if $w \cdot h = h$ for all $w \in W$, and it is $W$-anti-invariant if $w \cdot h = (-1)^{\sgn(w)}h$ for all $w \in W$. The set of $W$-invariants (resp. $W$-anti-invariants) form a subspace of $H$ that we denote $H^W$ (resp. $H^a$). A $W$-equivariant linear map $H \rightarrow H'$ induces linear maps $H^a \rightarrow (H')^a$ and $H^W \rightarrow (H')^W$.

\begin{lemma}\label{lem:exact}Taking invariants or anti-invariants is an exact functor from the category of $W$-vector spaces and $W$-equivariant morphisms to the category of vector spaces. 
\end{lemma}
\begin{proof}The only nontrivial part is to show that taking anti-invariants (resp. invariants) preserves surjections. For that, we use that the inclusions $H^a \to H$ (resp. $H^W \to H$) are split by the maps
\[
\phi^a(h) = |W|^{-1}\sum_{w \in W} sgn(w) w\cdot h,\qquad \phi^W(h) = |W|^{-1}\sum_{w \in W}w \cdot h.
\]
%In particular $h \in H$ is anti-invariant if and only if $\phi(h)=h$. 
which commute with $W$-equivariant linear maps. The  claim follows.
%which gives if $F: H \to H'$ is a linear $W$-equivariant map, then $F\circ \phi = \phi \circ F$. \Rachel{are y'all convinced? should we add more detail?}
%\Wendelin{The only nontrivial part is to show that taking anti-invariants (or invariants) preserves surjections. For this, consider any $W$-equivariant linear map $F \colon U \rightarrow V$. 
%The inclusion $V^a \to V$ is split by the map
%\[
%\phi(v) = |W|^{-1}\sum_{w \in W} sgn(w)v.
%\]
%and $\phi$ commutes with $F$. 
%In particular $h \in H$ is anti-invariant if and only if $\phi(h)=h$.}
\end{proof}

The application of this general discussion to the problem at hand is as follows. The $G$-action on $X$ induces a $W$-action on $X\sslash_\theta T$ commuting with the $K$-action. We therefore get an action of $W$ on $A^*_K(X\sslash_\theta T)_\QQ$ by defining $W$ to act on the classifying space $(X\sslash_\theta T) \times_K U_N$ by the action just described on the first factor and the trivial action on the second factor. 
% This is compatible with defining $W$ to act on $T \times K$ via its usual action on $T$ and the trivial action on $K$. \Rachel{what does that sentence mean, and where is the $W$-action used} Note that if $\cX = X\sslash_\theta T$ then the nonequivariant limit \eqref{eq:limit} is a $W$-equivariant ring homomorphism. \Rachel{is this used}
%Since $W$ acts trivially on $R_K$, the splittings of \eqref{eq:limit} generated by Lemma \ref{lem:lh} are also $W$-equivariant

\subsection{Equivariant Chow ring of $X\sslash_\theta G$}\label{sec:nonab chow ring}
We begin with a general observation: let $K$ and $T$ be tori and let $\cY$ be a smooth variety such that $K \times T$ acts on $\cY$.  Then for any integer $N$ we have
\begin{equation}\label{eq:chow1}
A^N_{K \times T}(\cY)_\QQ = A^N_K(\cY \times_T V_N)_\QQ
\end{equation} 
where $T$ acts by scaling  $V_N = (\mathbb{A}^{N+1}\setminus 0)^{\rk(T)}$ and $K$ acts trivially on $V_N$. If moreover the stack quotient $[\cY/T]$ is a scheme then the groups in \eqref{eq:chow1} are also equal to $A^N_K([\cY/T])_\QQ$.

Fix now a set $\Phi^+$ of positive roots of $G$ with respect to $T$ and set
\[
\Delta:= \prod_{\rho \in \Phi^+} c_1^K(\cL_\rho) \in A^*_{K \times T}(X)_\QQ %A^*_K(X\sslash_{\theta, G} T)_\QQ
\]
where $\cL_\rho$ is the $K$-linearized line bundle on $X\sslash_\theta T$ induced by the tensor product of the 1-dimensional representation of $T$ defined by $\rho$ with the trivial representation of $K$. Using \eqref{eq:chow1}, we also use $\Delta$ to denote the restriction of this class to $A^*_{K}(X\sslash_\theta T)_\QQ$ and $A^*_K(X\sslash_{\theta, G} T)_\QQ$. We have a diagram of morphisms of $R_K$-modules as follows.
\begin{equation}\label{eq:chow-diagram}
\begin{tikzcd}
& A^*_K(X\sslash_\theta T)_\QQ^W \arrow[r, "\cup \Delta"] \arrow[d, "j^*"] & A^*_K(X\sslash_\theta T)_\QQ^a \arrow[d, "j^*"]\\
A^*_K(X\sslash_\theta G)_\QQ \arrow[r, "g^*"] & A^*_K(X\sslash_{\theta, G} T)_\QQ^W \arrow[r, "\cup \Delta"] & A^*_K(X\sslash_{\theta, G} T)_\QQ^a
\end{tikzcd}
\end{equation}
The bottom horizontal arrows are isomorphisms by \cite[Prop~2.4.1, Lem~5.4.1]{nonab1}\footnote{The references a priori hold non-equivariantly, but the same proofs can be used to prove the results for $K$-equivariant Chow rings. One does this by replacing $Z^s(G)$ in the references by $X^s(G) \times_K U_N$. Note that $X^s(G)\times_K U_N$ is a scheme because $K$ is special (\cite[Prop~23]{EG98}). If we let $G$ act trivially on $U_N$, the induced $G$-scheme $X^s(G)\times_K U_N$ carries a $G$-linearized ample line bundle such that all points of $X^s(G)\times_K U_N$ are stable, so the hypotheses in \cite[(2.1)]{ES} are satisfied. This is all that is needed for the arguments in \cite[Prop~2.4.1, Lem~5.4.1]{nonab1} to go through.}.
Hence we may define two morphisms comparing the Chow groups of $X\sslash_\theta G$ and $X\sslash_\theta T$:
\begin{align*}
p^W:& \,\, A^*_K(X\sslash_\theta T)_\QQ^W \to A^*_K(X\sslash_\theta G)_\QQ & p^W &:= (g^*)^{-1} \circ j^* \\
p^a :& \,\, A^*_K(X\sslash_\theta T)_\QQ^a \to  A^*_K(X\sslash_\theta G)_\QQ & p^a &:= (g^*)^{-1} \circ (\cup \Delta)^{-1} \circ j^*
\end{align*}

\begin{proposition}
The morphism $p^W$ is a surjective $R_K$-algebra homomorphism with kernel equal to the kernel of $\cup \Delta$, and the morphism $p^a$ is an isomorphism of $R_K$-modules.
\end{proposition}

The proposition is an immediate consequence of the following lemma.

\begin{lemma}\label{lem:compare-chow}
In \eqref{eq:chow-diagram} we have that
\begin{enumerate}[(a)]
\item $\cup \Delta: A^*_K(X\sslash_\theta T)^W_\QQ \to A^*_K(X\sslash_\theta T)^a_{\QQ}$ is surjective, and
\item $j^*: A^*_K(X\sslash_\theta T)^a_\QQ \to A^*_K(X\sslash_{\theta, G} T)^a_\QQ$ is an isomorphism.
\end{enumerate}
\end{lemma}
\begin{proof}

For (a), the open embedding $X^{ss}_\theta(T) \to X$ induces a $W$-equivariant surjection from $A^*_{K\times T}(X)_\QQ$ to $A^*_{K\times  T}(X^{ss}_\theta(T))_\QQ$ and hence (using \eqref{eq:chow1}) a $W$-equivariant surjection
\[
Sym^\bullet(\chi(T)_\QQ) \otimes Sym^\bullet(\chi(K)_\QQ) \twoheadrightarrow A^*_K(X\sslash_\theta T)_\QQ.
\]
By Lemma \ref{lem:exact} we get a commuting square with surjective horizontal arrows
\[
\begin{tikzcd}
(Sym^\bullet(\chi(K)_\QQ) \otimes Sym^\bullet(\chi(T)_\QQ))^W \arrow[d, "\cup \Delta"] \arrow[r, twoheadrightarrow] & (A^*_K(X\sslash_\theta T)_\QQ)^W \arrow[d, "\cup \Delta"] \\
(Sym^\bullet(\chi(K)_\QQ) \otimes Sym^\bullet(\chi(T)_\QQ))^a \arrow[r, twoheadrightarrow] & (A^*_K(X\sslash_\theta T)_\QQ)^a
\end{tikzcd}
\]
Since $W$ acts trivially on $Sym^\bullet(\chi(K)_\QQ)$ and since $\Delta$ is contained in the subring $Sym^\bullet(\chi(T)_{\QQ})$ it is enough to show that
\[
\cup \Delta: (Sym^\bullet(\chi(T)_\QQ))^W \to (Sym^\bullet(\chi(T)_\QQ))^a
\]
is a surjection. This follows from \cite[Lem~1.2]{ES}.

For (b) we adapt the proof of \cite[Lem 3 on p.130]{brion91}. We begin by briefly recalling the Hesselink stratification,
%\Rachel{ask Vicky how to cite}, 
referring the reader to \cite[Sec 2]{hoskins} for the full definition. The Hesselink stratification is a decomposition
\begin{equation}\label{eq:strat}
X^{us}(G) = \bigsqcup_{\lambda \in \cB} \Sigma_{[\lambda]}
\end{equation}
where $\cB$ is a set of 1-parameter subgroups of $G$ and $\Sigma_{[\lambda]}$ is a locally closed $G$-invariant subvariety of $X$ depending only on the conjugacy class $[\lambda]$ of $\lambda$. There is a partial order on the index set in \eqref{eq:strat} such that if we add in $X^{ss}(G)$ as the lowest stratum, the unions $\bigcup_{[\lambda'] < [\lambda]} \Sigma_{[\lambda']}$ and $\bigcup_{[\lambda'] \leq [\lambda]} \Sigma_{[\lambda']}$ are open in $X$.

Moreover, the strata $\Sigma_{[\lambda]}$ are flag bundles as follows.  There is a closed subvariety $Y_{\lambda} \subset \Sigma_{[\lambda]}$ contained in the complement of $X^{ss}_\theta(T)$, as well as a parabolic subgroup $P_\lambda \subset G$ containing $T$, such that $Y_{\lambda}$ is $P_{\lambda}$-invariant and the map 
\begin{equation}\label{eq:strat-iso}
G \times_{P_{\lambda}} Y_{\lambda} \to \Sigma_{[\lambda]}
\end{equation}
induced by the $G$-action is a $T$-equivariant isomorphism.\footnote{A proof of the isomorphism \eqref{eq:strat-iso} in the projective setting can be found in \cite[Thm 13.5]{kirwan}. In our affine setting, one can modify the proof of \cite[Lem 13.4]{kirwan} to show that \eqref{eq:strat-iso} is a bijection on $\CC$-points, and then conclude using smoothness of $Y_\lambda$ and the action morphism $G \times X \to X$ combined with Zariski's Main Theorem (see \cite[\href{https://stacks.math.columbia.edu/tag/02LR}{Tag 02LR}]{stacks-project}).} Here we let $T$ act on the left hand side of \eqref{eq:strat-iso} via left multiplication on $G$.
In our situation, where we have torus $K$ whose action on $X$ commutes with the action of $G$, it is evident from the definitions (see e.g. \cite[Sec 2]{hoskins}) that the varieties $Y_\lambda$ and $\Sigma_{[\lambda]}$ are $K$-invariant; the isomorphism \eqref{eq:strat-iso} is clearly $K$-equivariant. We set 
\[
\Sigma^{ss}_{[\lambda]} = \Sigma_{[\lambda]} \cap X^{ss}(T).
\]

\begin{lemma}\label{lem:inductive-step}
We have $\Delta = 0$ in $A^*_K(\Sigma^{ss}_{[\lambda]}/T)$ for all $\lambda \in \cB$.
\end{lemma}
\begin{proof}
For simplicity set $\Sigma := \Sigma_{[\lambda]}$, $Y:= Y_\lambda$, and $P := P_\lambda$. There is a $T \times K$-equivariant commuting square
\[
\begin{tikzcd}
\Sigma^{ss} \arrow[r] \arrow[d] & G \times_P Y \arrow[d] \\
(G/P) \setminus (WP/P) \arrow[r] & G/P
\end{tikzcd}
\]
where $WP/P \subset G/P$ is the set of $T$-fixed points (equal to the quotient of the Weyl orbit of $P$). This is because $\Sigma^{ss}$ does not meet the subvariety $Y \simeq P \times_P Y \subset G \times_P Y$ and hence does not meet its $W$-orbit either. Since $A^*_K(\Sigma^{ss}/T)_\QQ = A^*_{K \times T}(\Sigma^{ss})_\QQ$ it is enough to show $\Delta=0$ in $A^*_{K \times T}((G/P)\setminus(WP/P))_\QQ$; hence, by excision, it is enough to show $\Delta$ is in the image of 
\[
i_*: A^*_T ((WP/P) \times_K U_N)_\QQ \to A^*_T( (G/P) \times_K U_N)_\QQ
\]
for sufficiently large $N$, where $i$ is the inclusion of the fixed locus.

By the localization formula in \cite[Thm 2]{EG98localization}, we may write $\Delta = i_*\alpha$ in $A^*_T( (G/P) \times_K U_N)_\QQ \otimes \Frac(R_T)$, for some class $\alpha \in A^*_T((WP/P) \times_K U_N)_\QQ.$ By Lemma \ref{lem:lh}, we see that $A^*_T( (G/P) \times_K U_N)_\QQ$ is a free $R_T$-module, so the natural homomorphism 
\[
A^*_T( (G/P) \times_K U_N)_\QQ \to A^*_T( (G/P) \times_K U_N)_\QQ \otimes \Frac(R_T)
\]
is injective.\footnote{
Some comments about why Lemma \ref{lem:lh} applies: $(G/P) \times_K U_N$ is isomorphic to the product of $G/P$ with a product of projective spaces. All factors have a $\Gm$-action with isolated fixed points. Moreover the Chow ring of $G/P$ has a $T$-equivariant basis of (Schubert) classes.
} We conclude that $\Delta = i_*\alpha$ in $A^*_T( (G/P) \times_K U_N)_\QQ.$
\end{proof}

To finish the proof of Lemma \ref{lem:compare-chow}, define $T\times K$-equivariant open subvarieties of $X^{ss}(T)$
\[
\Sigma^{ss}_{\leq [\lambda]} :=  \bigcup_{[\lambda'] \leq [\lambda] }\Sigma^{ss}_{[\lambda']} \quad \quad \text{and} \quad \quad \Sigma^{ss}_{< [\lambda]} :=  \bigcup_{[\lambda'] < [\lambda] }\Sigma^{ss}_{[\lambda']}.
\]
To show that the anti-invariant part of $j^*$ is an isomorphism, first note that it is a surjection because it is given by restriction to an open subset. To prove injectivity, note that both morphisms in the composite
\begin{equation}\label{eq:p9}
\begin{tikzcd}
(Sym^\bullet(\chi(T)_\QQ) \otimes Sym^\bullet(\chi(K)_\QQ))^a \arrow[r] & A^*_K(X\sslash_\theta T)_\QQ^a \arrow[r, "j^*"] & A^*_K(X\sslash_{\theta, G} T)_\QQ^a
\end{tikzcd}
\end{equation}
are surjective, as each is induced by restriction to an open substack.
Moreover, for each $[\lambda]$ appearing in \eqref{eq:strat} the morphism $j^*$ in \eqref{eq:p9} further factors as a composition of morphisms
\[
A^*_K(X\sslash_\theta T)_\QQ^a \to A^*_K( \Sigma^{ss}_{\leq [\lambda]} / T)_\QQ^a \to A^*_K(\Sigma^{ss}_{<[\lambda]}/T)_\QQ^a \to A^*_K(X\sslash_{\theta, G} T)_\QQ^a
\]
each of which is again given by restriction to an open subvariety and therefore is surjective.
%In particular the morphism in part (b) of Lemma \ref{lem:compare-chow} is surjective. 
To show that the anti-invariant part of $j^*$ is injective, suppose 
\begin{equation}\label{eq:compare-chow10}
\alpha \in (Sym^\bullet(\chi(T)_\QQ) \otimes Sym^\bullet(\chi(K)_\QQ))^a
\end{equation}
has the property that its image in $A^*_K(X\sslash_{\theta, G} T)_\QQ^a$ is zero. It is enough to show that its image in $A^*_K(X\sslash_\theta T)_\QQ$ is also zero. We argue by induction: if the image of $\alpha$ in $A^*_K(\Sigma^{ss}_{< [\lambda]}/T)_\QQ^a$ is zero, then by excision we see that $\alpha \in A^*_K(\Sigma^{ss}_{[\lambda]}/T)_\QQ^a$. Since $\alpha$ is an anti-invariant of \eqref{eq:compare-chow10}, \cite[Lem~1.2]{ES} tells us that $\Delta$ divides $\alpha$, so Lemma \ref{lem:inductive-step} implies $\alpha=0$ in $A^*_K(\Sigma^{ss}_{[\lambda]}/T)_\QQ^a.$ Hence $\alpha=0$ in $A^*_K(\Sigma^{ss}_{\leq [\lambda]}/T)_\QQ^a$ as well.
\end{proof}

\section{Abelianization of quasimap $I$-functions}\label{sec:abelianization}

In this section we reformulate the abelianization theorem in \cite[Cor 6.3.1]{nonab1} to parallel the statements in \cite[Thm 6.1.2]{frobman} and \cite[p. 103]{twoproofs}.\footnote{We can almost just cite \cite[Thm 6.1.2]{frobman} for our use here, except that this reference works only with compact targets.}
Let $(X, G, \theta, K)$ be as in Section \ref{sec:givental}. Let $T \subset G$ be a maximal torus and let $W$ be the Weyl group of $T$ in $G$. For $\tilde \beta \in \Hom (\chi(T), \ZZ)$ define
\[
I_{\tilde \beta} := \prod_{\xi} \frac{\prod_{\ell=-\infty}^0 (c_1^K(\cL_\xi) + \ell z)}{\prod_{\ell=-\infty}^{\tilde \beta(\xi)} (c_1^K(\cL_\xi) + \ell z)} \quad \quad \in A^*_K(X\sslash_\theta G)_\QQ((z^{-1}))
\]
where the product runs over all weights $\xi$ of the $T \times K$ action on $X$, we denote by $c_1^K$ the equivariant first Chern class, and $\cL_\xi$ is the $K$-equivariant line bundle on $X\sslash_\theta G$ induced by $\xi$.

Let $k$ denote the rank of $T$, fix $pr_1, \ldots, pr_k$ an integral basis for $\chi(T)$, and set %\Qaasim{I've gotten a bit confused about why these are supposed to be non-equivariant cohomology classes, even though the whole I-function takes values in the equivariant cohomology}\Rachel{better?} 
$H_i := c_1^K(\cL_{pr_i})$, where $\cL_{pr_i}$ is the $K$-equivariant line bundle on $X\sslash_\theta G$ induced by $pr_i$ and the trivial character of $K$. Let $\log(y_1), \ldots, \log(y_k)$ be formal variables and define
\[
\square := \sum_{i=1}^k H_i \log(y_i). 
\]
The Weyl group acts by group automorphisms on $\chi(T)$ and hence on the classes $H_i$. Let $\{P_j(H_1, \ldots, H_k)\}_{j \in J}$ be a set of $W$-invariant polynomials in the $H_i$ and let $\{\bx_j\}_{j \in J}$ be corresponding formal variables. For effective $\tilde \beta \in \Hom(\chi(T), \ZZ)$ we let $q^{\tilde \beta}$ denote the corresponding element of $\QQ[\Hom(\chi(T), \ZZ)]$. Define a series 
\begin{equation}\label{eq:nonab1}
\II_T := z e^{\square/z} \sum_{\tilde \beta \in \Hom(\chi(T), \ZZ)} q^{\tilde \beta} e^{\sum_{i=1}^k \log(y_i) \tilde \beta(pr_i)} e^{\sum_{j \in J} \bx_j P_j(\bH + \tilde \beta z)/z} I_{\tilde \beta}
\end{equation}
where we set 
\[
P_j(\bH + \tilde \beta z) = P_j(H_1 + \tilde \beta(pr_1)z, \ldots, H_k+\tilde \beta(pr_k)z).
\]

%\Rachel{added this paragraph. does it look ok?}\Wendelin{Yes}
Recall that the Novikov ring $\Lambda_{F_K}$ is a completion of
$F_k[\mathrm{Eff}(X\sslash_\theta G)]
\subseteq F_K[\mathrm{Hom}(\chi(G),
\mathbb{Z})]$; for effective $\beta \in \Hom(\chi(G), \ZZ)$ let $q^\beta$ denote the corresponding element of $\Lambda_{F_K}$, and let $q^{\tilde \beta} \mapsto q^\beta$ denote the specialization induced by $\Hom(\chi(T), \ZZ) \to \Hom(\chi(G), \ZZ)$. Let $\Phi_+$ be a set of positive roots, let $\zeta$ be their sum, and let $\Delta$ be as in Section \ref{sec:nonab chow ring}. Finally, for any $\rho \in \chi(T)$ let $a_i^\rho$ be integers such that $\rho = \sum_{i=1}^k a^\rho_i pr_i$ and set
\[
\partial_\Delta = \prod_{\rho \in \Phi_+} \partial_\rho \quad \quad \quad \quad\text{where} \quad \partial_\rho = z \sum_{i=1}^k a_i^\rho \partial_{\log(y_i)}.
\]

The two formulas in the following theorem are inspired by the two distinct formulas in \cite[Thm 6.1.2]{frobman} and \cite[p. 103]{twoproofs}, respectively.

\begin{theorem}[Nonabelian big $I$-functions, {\cite[Thm 3.3]{bigI} \cite[Cor 6.3.1]{nonab1}}]\label{thm:nonabelianI}
Let $y \mapsto \by$ be a linear specialization of variables with the property that $\square|_{y \mapsto \by}$ is Weyl-invariant (where $W$ acts on the Chow group). 
%Let $\zeta = \sum_{\rho \in \Phi^+}{ \rho}$.
%and write $\zeta = \sum_{i=1}^k a^{\zeta}_i pr_i$ for some unique integers $a^{\zeta}_i$.
Then there is a unique series $\II_G$ with the property that 
\begin{equation}\label{eq:def IG}
\Delta g^*\II_G =  \partial_\Delta j^*\II_T|_{y=\by, \; q^{\tilde \beta} = (-1)^{\tilde \beta(\zeta)}q^\beta}.
\end{equation}
If $a_i \in \ZZ$ are the unique integers satisfying $\zeta = \sum_{i=1}^k a_i pr_i$ and if $\log(y)= \log(\by) + i \pi \ba$ denotes the change of variables $\log(y_i) \mapsto \log(y_i) + i\pi a_i$ followed by $y \mapsto \by$, then \eqref{eq:def IG} is also equal to 
\begin{equation}\label{eq:def IG2}
e^{-i\pi c_1^K(\cL_\zeta)}\partial_\Delta j^*\II_T|_{\log(y)=\log(\by)+ i \pi \ba, \; q^{\tilde \beta} = q^\beta}
\end{equation}
when $\cL_\zeta$ is equipped with the $K$-linearization induced by the trivial character of $K$. Moreover, if the $K$-action on $X\sslash_\theta G$ has isolated fixed points and isolated 1-dimensional orbits then $\II_G$ is on the cohomological Lagrangian cone of $X\sslash_\theta G$.
\end{theorem}
The ``cohomological Lagrangian cone'' referred to in Theorem \ref{thm:nonabelianI} is the cone defined using $H^*_K(X\sslash_\theta G; \QQ)$ rather than $A^*_K(X\sslash_\theta G)_\QQ$, and when we say $\II_G$ is an element we really mean that if we apply the Borel-Moore homology valued cycle map and then Poincar\'e duality isomorphism to $\II_G$, the image is on the cone. We expect Theorem \ref{thm:nonabelianI} holds true when ``cohomologial Lagrangian cone'' is replaced with ``Lagrangian cone.'' The following remark further clarifies the statement of Theorem \ref{thm:nonabelianI} and will also be used in its proof.
\begin{remark}\label{rmk:choice of roots}
The series $\II_G$ is independent of the choice of $\Phi_+$, which a priori affects the definition of $\Delta$, $\partial_\Delta$ and $\zeta$. To see this, suppose $\Phi_+'$ is another set of positive roots, and let $\Delta', \partial_{\Delta'}, \zeta'$ be the corresponding objects for $\Phi_+'$. If $R$ is the set of $\rho \in \Phi_+$ such that $-\rho$ is in $\Phi_+'$, then we have
\begin{equation}\label{eq:change sign}
\Delta = (-1)^{|R|}\Delta' \quad \quad \quad \partial_{\Delta} = (-1)^{|R|}\partial_{\Delta'} \quad \quad \quad \zeta = \zeta' + 2\sum_{\rho \in R} \rho.\end{equation}
It follows that when we replace $\Delta$ by $\Delta'$ and $\partial_{\Delta}$ by $\partial_{\Delta'}$, the left and right hand sides of \eqref{eq:def IG} change by the same sign. On the other hand, we see that the right hand side of \eqref{eq:def IG} does not change when we replace $\zeta$ by $\zeta'$. (In fact, we may use any $\zeta''\in \chi(T)$ in place of $\zeta$ as long as $\zeta''$ satisfies $\tilde \beta(\zeta) - \tilde \beta(\zeta'') \in 2\ZZ$ for all effective $\tilde \beta$.)

\end{remark}

\begin{proof}[Proof of Theorem \ref{thm:nonabelianI}]

We show that $\II_G$ agrees with the (Chow version of the) series in \cite[Cor~6.3.1]{nonab1}, which for this proof we denote $\mathbf{I}$; the cohomological version of this series is on the Lagrangian cone when the $K$-action has isolated fixed points and 1-dimensional orbits by \cite[Thm 3.3]{bigI}. By the definition in \cite[Cor~6.3.1]{nonab1}, we have
\begin{equation}\label{eq:def other I}
\begin{gathered}
g^*\bI = \left( z \sum_{\tilde \beta} q^{\tilde \beta} e^{\sigma_{\tilde \beta}/z} \prod_\rho \frac{\prod_{\ell=-\infty}^{\tilde \beta(\rho)} (c_1^K(\cL_\rho)  + \ell z)}{\prod_{\ell=-\infty}^{0} (c_1^K(\cL_\rho)  + \ell z)} j^*I_{\tilde \beta} \right)|_{y \mapsto \by, \;q^{\tilde \beta} \mapsto q^\beta}\\
\text{where}\quad \sigma_{\tilde \beta} := \sum_{i=1}^k \log(y_i)(H_i + \tilde \beta(pr_i)z) + \sum_{j \in J} \bx_jP_j(\bH+\tilde \beta z)
\end{gathered}
\end{equation}
and the product ranges over all roots $\rho$ of $G$ with respect to $T$, and $\cL_\rho$ has the $K$-linearization induced by the trivial character of $K$.
We compute the abelianization factor:
\begin{equation}\label{eq:factor}
 \prod_\rho \frac{\prod^{\tilde \beta(\rho)}_{\ell=-\infty}(c_1^K(\cL_{\rho})  + \ell z)}{\prod^{0}_{\ell=-\infty}(c_1^K(\cL_{\rho})  + \ell z)} = \prod_{\tilde \beta(\rho) > 0} \left((-1)^{\tilde \beta(\rho)}\frac{c_1^K(\cL_{\rho})+ \tilde \beta(\rho)z}{c_1^K(\cL_\rho)}\right).
\end{equation}
Note that we may add any number of factors with $\tilde \beta(\rho)=0$ to the product on the right hand side without changing it, since if $\tilde \beta(\rho)=0$ the corresponding factor is 1. But there is a set of positive roots $\Phi_+^{\tilde \beta}$ containing $\{ \rho \mid \tilde \beta(\rho) > 0\}$ and contained in $\{ \rho \mid \tilde \beta(\rho) \geq 0\}$%(\Wendelin{Is the point that $\tilde{\beta})$ definers a halfspace in $\chi(T)$ and we can always choose a set of positive roots contained in a given halfspace?)}\Rachel{yes. we're also using that such a set of positive roots will always contain the roots that lie strictly on the positive side of the hyperplane}
, so we may take the above product over $\Phi_+^{\tilde \beta}$. If we let $\Delta_{\tilde \beta}$ denote the fundamental Weyl anti-invariant class corresponding to $\Phi_+^{\tilde \beta}$ then as in Remark \ref{rmk:choice of roots} we have $\Delta = {\epsilon(\tilde \beta)}\Delta_{\tilde \beta}$ for some $\epsilon(\tilde \beta) \in \{-1, 1\}$. Hence \eqref{eq:def other I} is equal to
\[
\Delta^{-1}\left(z \sum_{\tilde \beta}q^{\tilde \beta} e^{\sigma_{\tilde \beta}/z} \epsilon(\tilde \beta) \prod_{\rho \in \Phi^{\tilde \beta}_+} (-1)^{\tilde \beta(\rho)} (c_1^K(\cL_{\rho})+ \tilde \beta(\rho)z) j^*I_{\tilde \beta}|_{y \mapsto \by, \;q^{\tilde \beta} \mapsto q^\beta} \right)
\]
We apply Remark \ref{rmk:choice of roots} two more times to write
\[\prod_{\rho \in \Phi^{\tilde \beta}_+}(-1)^{\tilde \beta(\rho)}=(-1)^{\tilde \beta(\zeta)} \quad \quad \text{and}\quad \quad \epsilon(\tilde \beta) \partial_{\Delta_{\tilde \beta}} = \partial_{\Delta}\]
so that \eqref{eq:def other I} now becomes
\begin{equation}\label{eq:hey}
\Delta^{-1}\left(z \sum_{\tilde \beta}q^{\tilde \beta} e^{\sigma_{\tilde \beta}/z}  \prod_{\rho \in \Phi_+} (-1)^{\tilde \beta(\rho)} (c_1^K(\cL_{\rho})+ \tilde \beta(\rho)z) j^*I_{\tilde \beta}|_{y \mapsto \by, \;q^{\tilde \beta} \mapsto q^\beta} \right).
\end{equation}
Here, given a Weyl
anti-invariant class
$x$, we write
$\Delta^{-1}x$ for the
preimage of $x$
under the
isomorphism $\cup
\Delta:
A^*_K(X\sslash_{\theta, G} T)^W_\QQ \to
A^*_K(X\sslash_{\theta, G} T)^a_\QQ$ in \eqref{eq:chow-diagram}. We now
show that the quantity
in parentheses in
\eqref{eq:hey} is anti-invariant; since this
quantity is just the
right hand side of
\eqref{eq:def IG}, this
will also prove the first
statement of the
theorem. We will prove this anti-invariance next; since the quantity in parentheses is just the right hand side of \eqref{eq:def IG}, this will also prove the first statement of the theorem.

For the claimed anti-invariance, since $W$ acts on the set of $\tilde \beta$ mapping to $\beta$ 
it is enough to show that the quantities
\[
B_{\tilde \beta}:= e^{\sigma_{\tilde \beta}/z}\prod_{\rho \in \Phi_+} (-1)^{\tilde \beta(\rho)} (c_1^K(\cL_\rho) + \tilde \beta(\rho)z) j^*I_{\tilde \beta}|_{y \to \by}
\]
have the property that whenever $w \in W$ is a reflection along a root, we have 
\begin{equation}\label{eq:wB}
w \cdot B_{\tilde \beta} = -B_{w\tilde \beta}.\end{equation}
Here we are using the natural action of $W$ on $\Hom(\chi(T), \ZZ)$ by the rule $(w\tilde \beta)(\xi) = \tilde \beta(w^{-1}\xi)$ for $\xi \in \chi(T)$.

The proof of \cite[Lem 5.4.2]{nonab1} shows that $w \cdot I_{\tilde \beta} = I_{w\tilde \beta}$. We next show that $w \cdot \sigma_{\tilde \beta} = \sigma_{w \cdot \tilde \beta}$. It is enough to show that for any $W$-invariant polynomial $P(H_1, \ldots, H_k)$ (where $W$ acts on the $H_i$) we have
\begin{equation}\label{eq:wP}
w\cdot P(\bH + \tilde \beta z) = P(\bH + w\tilde \beta z)
\end{equation}
where on the left hand side, $w$ acts only on elements of Chow. There is a $k\times k$ integer matrix $A_w = (a_{ij})$ with the property that $wH_j = \sum_{i=1}^k a_{ij}H_i$. If we view $\bH + \tilde \beta z$ as the vector $(H_1 + \tilde \beta(pr_1)z, \ldots, H_k + \tilde \beta(pr_k)z)$, then Weyl invariance of $P$ tells us
\[
P(\bH + \tilde \beta z) = P(A_w(\bH + \tilde \beta z)) = P(A_w\bH + A_w \tilde \beta z).
\]
But it is easy to see that $A_w\tilde \beta$ is the vector whose entries are $(w^{-1}\tilde \beta(pr_1), \ldots, w^{-1}\tilde \beta(pr_k))$. It follows that \eqref{eq:wP} holds.

Finally, to complete the proof of \eqref{eq:wB} we compute
\begin{align*}
w \cdot \prod_{\rho \in \Phi_+} (-1)^{\tilde \beta(\rho)} (c_1^K(\cL_\rho) + \tilde \beta(\rho)z)  &= \prod_{\rho \in \Phi_+} (-1)^{\tilde \beta(\rho)} (c_1^K(\cL_{w\rho}) + \tilde \beta(\rho)z) \\
&= \prod_{\rho \in w\Phi_+} (-1)^{w\tilde \beta(\rho)} (c_1^K(\cL_\rho) + w\tilde \beta(\rho)z).
%&= \prod_{\rho \in w\Phi^+ \cap \Phi^+} (-1)^{w\tilde \beta(\rho)} (c_1(\cL_\rho) + w\tilde \beta(\rho)z)\prod_{\rho \in w\Phi^+ \cap \Phi^-} (-1)^{w\tilde \beta(\rho)} (c_1(\cL_\rho) + w\tilde \beta(\rho)z)
\end{align*}
We break this product up into two, one over $\rho \in w\Phi_+ \cap \Phi_+$ and one over $\rho \in w\Phi_+ \cap \Phi_-$, where $\Phi_- := -\Phi_+$. We then reindex the second factor by replacing $\rho$ with $-\rho$ and find that this product is equal to
\[
\prod_{\rho \in w\Phi_+ \cap \Phi_+} (-1)^{w\tilde \beta(\rho)} (c_1^K(\cL_\rho) + w\tilde \beta(\rho)z)\prod_{\rho \in -w\Phi_+ \cap \Phi_+} (-1)^{-w\tilde \beta(\rho)} (-c_1^K(\cL_\rho) - w\tilde \beta(\rho)z) \]
\[= (-1)^{\#( -w \Phi_+ \cap \Phi_+)}\prod_{\rho \in \Phi_+} (-1)^{w\tilde \beta(\rho)} (c_1^K(\cL_\rho) + w\tilde \beta(\rho)z)
\]
where in the equality we used that $-w\Phi_+ = w\cdot(-\Phi_+) = w\Phi_-$, so that $w\Phi_+$ and $-w\Phi_+$ partition the set of roots. 
%(\Wendelin{Bit confused here: it looks to me that there are two minus signs in the second term which should cancel out?}\Rachel{not sure what the confusion is. ask in the next mtg?}) 
It remains to show that the cardinality of $ -w\Phi_+ \cap \Phi_+$ is odd. This is true because $-w$ is an involution on this set with a unique fixed point (the simple root along which $w$ is a reflection).
From here, we note that $\Delta$ times \eqref{eq:hey} is equal to the right hand side of \eqref{eq:def IG} and to \eqref{eq:def IG2} by direct computation.

\end{proof}

\section{The Grassmannian flop}\label{sec:prooftop}
In this section we let $X, G, \theta_{\pm}$ and $K$ be as in Section \ref{sec:statement} and we prove Theorem \ref{thm:main}. 
That is, $X=M_{k \times n} \times M_{n \times k}$ and $G = GL(k)$ acts on $X$ by 
\[
g\cdot(x,y) = (gx, yg^{-1})
\]
for $g \in G$ and $(x, y) \in M_{k\times n} \times M_{n \times k}.$
%The integral characters of $G$ are powers of the determinant character, and the tensor product of this lattice with $\QQ$ is a 1-dimensional vector space. If we 
The characters $\theta_{\pm}$ of $G$ are given by 
\[\theta_+(g) = \det(g) \quad \quad \quad \quad \theta_-(g) = \det(g)^{-1},\] 
and the torus $K = (\CC^*)^{2n}$ acts on $X$ by
\[
(k_1, k_2)\cdot(x,y) = (xk_1^{-1}, k_2y) \quad \quad \quad \quad (k_1, k_2) \in (\CC^*)^n \times (\CC^*)^n.
\]

Let $T \subset G = GL(k)$ be the subgroup of diagonal matrices and let $W=S_k$ be the Weyl group. Our argument consists of the following five steps, each of which will be addressed in its own section.
\begin{enumerate}
\item Define a toric variety $\mathcal{M}_T$ and an action of $W$ on $\mathcal{M}_T$, and define $\mathcal{M}_G$ to be the $W$-invariant subvariety. (The variety $\cM_T$ is a compactification of the stringy K\"ahler moduli space and arises naturally in the theory of variation of GIT as the secondary toric variety of $X\sslash_{\theta_\pm} T$). Define certain torus fixed points $p_0, \ldots, p_k \in \cM_T$ corresponding to characters $\theta_0, \ldots, \theta_k$ of $T$, such that the points $p_0, p_k$ are in $\cM_G$ and $\theta_0, \theta_k$ are restrictions of characters of $G$. The GIT quotients $X\sslash_{\theta_0} G$ and $X\sslash_{\theta_k} G$ will be the two sides of the Grassmannian flop.
\item Define a ``small $I$-function'' $I^c_T$ for $c=0, \ldots, k$ such that each $I^c_T$ is on the Lagrangian cone of $X\sslash_{\theta_c} T$ and defines a multi-valued analytic function on a punctured neighborhood of $p_c \in \cM_T$. 
%(\Wendelin{small point: should we say multivalued on a punctured neighbourhood?}
\item Show that the series $I^0_T$ and $I^k_T$ are related by analytic continuation on $\cM_T$, and that the path of analytic continuation can be homotoped to lie in $\cM_G$.
\item Define ``big $I$-functions'' $\II^0_T$ and $\II^k_T$ that are related to $I^0_T$ and $I^k_T$ via explicit reconstruction and show that these are related by the same analytic continuation as $I^0_T$ and $I^k_T$.
\item Define series $\II^-_G$ and $\II^+_G$ that are related to $\II^0_T$ and $\II^k_T$ via abelianization and show that these are related by analytic continuation.
\end{enumerate}

In step 3 we reduce to the analytic continuation in \cite{coates2018crepant}. This reduction is the most technical part of the paper.

\subsection{Step 1: Variation of GIT}\label{sec:step1}

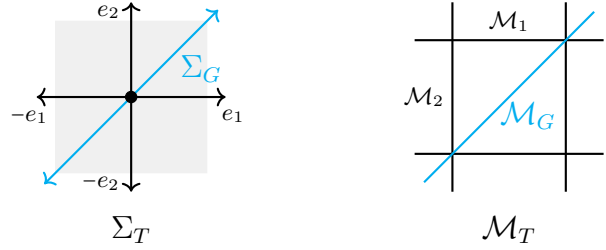
\begin{figure}[ht]
    \centering
    \[
    \begin{tikzpicture}[scale=.25]
    \node at (0, -7) {\Large $\Sigma_T$};
    \draw [gray!12, fill] (-4, -4) rectangle (4, 4);
\draw [thick, <->] (-5, 0)--(5, 0);
\draw [thick, <->]  (0, -5)--(0, 5);
\draw [cyan, thick, <->]  (-4.6, -4.6)--(4.6, 4.6);
\node [cyan] at (3.7,1.5) {\Large $\Sigma_G$};
\draw [fill] (0,0) circle [radius=.3];
\node at (5.4, -1) {\small$e_1$};
\node at (-1.6, -4.5) {\small$-e_2$};
\node at (-1.2, 4.5) {\small$e_2$};
\node at (-5.4, -1) {\small$-e_1$};

\begin{scope}[shift={(20,0)}]
\node at (0, -7) {\Large $\cM_T$};
\draw [thick] (-5, 3)--(5, 3);
\draw [thick] (-5, -3)--(5, -3);
\draw [thick] (-3, -5)--(-3, 5);
\draw [thick] (3, -5)--(3, 5);
\draw [thick, cyan] (-4.5, -4.5)--(4.5, 4.5);
\node at (0, 4) {\small$\cM_1$};
\node at (-4.5, 0) {\small$\cM_2$};
\node [cyan] at (1, -1) {\Large $\cM_G$};
\end{scope}
    \end{tikzpicture}
    \]
    \caption{The left diagram depicts the fan $\Sigma_T \subset \chi(T)_{\RR}$ when $k=2$ (consisting of four 2-dimensional cones and their faces). The Weyl group $S_2$ acts by reflection over the diagonal, and $\Sigma_G$ is the fan pictured with two 1-dimensional cones and the origin. The right diagram is the toric variety $\cM_T \simeq \PP^1 \times \PP^1$ with its four torus-invariant divisors and the Weyl-invariant subspace $\cM_{G} \simeq \PP^1$. The divisors labeled $\cM_1$ and $\cM_2$ will appear in the analytic continuation in Section \ref{sec:cij}. }
    \label{fig:only}
\end{figure}

Let $\chi(T)$ denote the character lattice of $T$, let $\chi(T)_\QQ$ be the associated $\QQ$-vector space, and define $\chi(G)$ and $\chi(G)_\QQ$ similarly. Our identification of $T$ with diagonal matrices in $GL(k)$ defines a basis $e_1, \ldots e_k$ of projection characters for $\chi(T)_\RR$. Let $\Sigma_T$ be the fan in $\chi(T)_\RR$ defined by the variation of GIT quotients corresponding to the pair $(X, T)$; that is, if $\theta$ is any character of $T$, then $\Sigma_T$ is the secondary fan of the toric variety $X\sslash_\theta T$. Explicitly, the rays of $\Sigma_T$ are $\{\pm e_i\}_{i=1}^k$ and the $d$-dimensional cones are the cones generated by $d$ linearly independent rays. (See Figure \ref{fig:only} for the case $k=2$.) 
%\Wendelin{In the left picture, the arrow pointing upwards should be labeled $e_2$!}
The Weyl group $W = S_k$ acts on $\chi(T)_\RR$ and on $\Sigma_T$. The Weyl-invariant subspace of $\chi(T)_\RR$ is the diagonal, which we may identify with $\chi(G)_\RR$, and the intersection of this subspace with $\Sigma_T$ defines a fan $\Sigma_G$ in $\chi(G)_\QQ \simeq \RR$ which is just the fan of $\PP^1$. Note that $\Sigma_G$ also gives the wall-and-chamber decomposition of $\chi(G)_\RR$ corresponding to the various possible GIT quotients.
We define $\cM_T$ and $\cM_G$ to be the toric varieties associated to the fans $\Sigma_T$ and $\Sigma_G$, respectively; that is,
\[
\cM_T = (\PP^1)^k \quad \quad \quad \quad \cM_G = \PP^1.
\]
Note that $W$ acts on $\cM_T$ and that $\cM_G$ is the Weyl-invariant subvariety.

The full-dimensional cones of $\Sigma_T$ are in bijection with GIT subvarieties $X\sslash_\theta T \subset [X/T]$ and also with torus fixed points in $\cM_T$. We will use $k+1$ of these cones, identified by characters $\theta_c$ for $c=0, \ldots, k$ where
\[
\theta_c := -e_1 - \ldots - e_c + e_{c+1} + \ldots + e_k.
\]
We let $p_c$ be the corresponding fixed point in $\cM_T$. Let $\theta_- = \theta_k$ and $\theta_+ = \theta_0$. Observe that $p_0$ and $p_k$ are contained in $\cM_G$ and that $\theta_+$ (resp. $\theta_-$) is the restriction of the determinant character of $G$ (resp. inverse of the determinant character of $G$).

We introduce coordinates on $\cM_T$ and $\cM_G$ as follows. Let $y_i^+$ be a variable corresponding to the dual of $e_i \in \chi(T)$ and let $y_i^-$ corrrespond to the dual of $-e_i$. The toric variety $\cM_T = (\PP^1)^k$ contains open sets
\[
U_c := \Spec(\CC[y_1^{-}, \ldots, y_c^-, y_{c+1}^+, \ldots, y_k^{+}])
\]
for $c=0, \ldots, k$. Note that $U_c$ contains $p_c$ as the origin in affine space. We abbreviate the coordinates on $U_c$ by writing $U_c = \Spec(\CC[y_i^{\epsilon^c_i}])$, where $\epsilon^c_i = -$ for $i=1, \ldots, c$ and $\epsilon^c_i=+$ for $i=c+1, \ldots, k$. Likewise, let $\by^+$ be a variable corresponding to the dual of the determinant character in $\chi(G)$ and let $\by^-$ correspond to the dual of the inverse of the determinant. The toric variety $\cM_G =\PP^1$ is covered by the open sets $U_{\pm} := \Spec(\CC[\by^{\pm}])$. 

In what follows we will regard $\cM_T$ and $\cM_G$ as complex analytic varieties, writing $\cO_{\cM_T}$ and $\cO_{\cM_G}$ for the corresponding sheaves of analytic functions. For each $c=0, \ldots, k$ we have the Givental module for $X\sslash_{\theta_c} T$ as well as an associated sheaf on $\cM_T$:
\[
H^c_T := F_K((z^{-1})) \otimes_{R_K} A^*_K(X\sslash_{\theta_c} T)_\QQ \quad \quad \quad \quad \cH^c_T := \cO_{\cM_T} \otimes_\CC H^c_T.
\]
If we fix a basis $\{\phi_i\}_{i=1}^{n^k}$ for $A^*(X\sslash_{\theta_c} T)_\QQ$ and equivariant parameters $\{\sigma_j\}_{j=1}^{2k}$ for $K$, then sections of $\cH^c_T$ on an open set $V \subset U_{c}$ may be written as formal series in the $\phi_i$, $\sigma_j$, $z^{-1}$ and $y^{-}_{1}, \ldots, y^-_d, y^+_{d+1}, \ldots, y^+_k$ with the property that the coefficient of a monomial in $\phi_i$, $\sigma_j$, $z^{-1}$ is a convergent power series in the $y^{\pm}_{i}$ on the domain $V$.
Likewise, we have the Givental modules for $X\sslash_{\theta_{\pm}} G$ as well as associated sheaves on $\cM_G$:
\[
H^{\pm}_G := F_K((z^{-1})) \otimes_{R_K} A^*_K(X\sslash_{\theta_{\pm}} G)_\QQ \quad \quad \quad \quad \cH^\pm_G := \cO_{\cM_G} \otimes_\CC H^\pm_G.
\]

\subsection{Step 2: Small toric $I$-functions}

For any index $c=0, \ldots, k$, let $H_i$ denote the element of $A^*_K(X\sslash_{\theta_c} T)_\QQ$ equal to the equivariant first Chern class of the line bundle induced by the projection character $e_i$ of $T$ and the trivial character of $K$.
Moreover let
\[R_K = \QQ[\lambda_1, \ldots, \lambda_n, \sigma_1, \ldots, \sigma_n]\]
denote the equivariant cohomology ring of a point, where $\lambda_1, \ldots, \lambda_n$ (resp. $\sigma_1, \ldots, \sigma_n$) are the first Chern classes for the line bundles on the classifying space $BK$ corresponding to the first (resp. last) $n$ projection characters of $K$.

 Introduce the Novikov ring $\Lambda^c_{F_K} = F_K\llbracket q_1^{-1}, \ldots,q_c^{-1}, q_{c+1},\ldots, q_k\rrbracket.$
 Finally, given integers $d_1, \ldots, d_k$, let
$\bd := (d_1, \ldots, d_k).$
  We then define 
\begin{align*}
I^c_T := z e^{\square_c/z}& \sum_{\substack{d_1, \ldots, d_c \leq 0 \\ d_{c+1}, \ldots, d_k \geq 0 }} 
\prod_{i=1}^c(y_i^-)^{-d_i} \prod_{i=c+1}^k(y_i^+)^{d_i}
\prod_{i=1}^k q_i^{d_i} I^c_{T, \bd}\\
\text{where}  \quad \quad \quad \square_c &= \sum_{i=1}^c \log(y_i^-)(-H_i) + \sum_{i=c+1}^k \log(y_i^+)(H_i)\\
\text{and} \quad \quad \quad I^c_{T, \bd} &= \prod_{j=1}^{n} \left(\prod_{i=1}^c\prod_{\ell=d_i+1}^0 (H_{i} -\lambda_j+\ell z)
\prod_{i=c+1}^{k} \prod_{\ell=1}^{d_i}(H_i-\lambda_j+\ell z)^{-1} \right)\\
& \quad \prod_{j=1}^n\left(\prod_{i=1}^c \prod_{\ell=1}^{-d_i}(-H_i +\sigma_j+ \ell z)^{-1}
\prod_{i=c+1}^k\prod_{\ell=-d_i+1}^{0}(-H_i +\sigma_j+ \ell z)\right).
\end{align*}

\begin{remark}\label{rmk:agrees}
The cycle map $A^*_K(X \sslash_{\theta_c} T)_\QQ \to H^*_K(X\sslash_{\theta_c} T; \QQ)$ is an isomorphism, as can be seen from the presentation of $H^*_K(X\sslash_{\theta_c}T; \QQ)$ in \cite[Sec 4.3]{coates2018crepant} and a direct computation of $A^*_K(X\sslash_{\theta_c} T)_\QQ$. (To compute the Chow group, one can use the exact sequence coming from the open embedding $X^{ss}_{\theta_c}(T) \subset X$.) 

After making this identification of Chow and cohomology, the series $I^c_T$ agrees with the one in \cite[(5.11)]{coates2018crepant}, up to a constant in the equivariant parameters, when one uses the GIT data $(X, T, \theta_c)$ and when one chooses the basis 
\[-e_1, \ldots, -e_c, e_{c+1}, \ldots, e_{k}\] for $\chi(T)$ (this basis is the one denoted in \cite{coates2018crepant} by $\mathrm{p}_1^+, \ldots, \mathrm{p}_r^+$).\footnote{Note that the ordering of our basis does not agree with the ordering of the set $\{\mathrm{p}_i^+\}_{i=1}^r$ in \cite[1038]{coates2018crepant}. This difference does not affect the $I$-function formula.} 
\end{remark}
We will interpret the series $I^c_T$ in two ways: formally algebraically, as a formal series in $\log(y_i^{\pm}), z^{\pm 1}$, and $q_i^{\pm 1}$ with coefficients in $F_K \otimes_{R_K} A^*_K(X\sslash_{\theta_c} T)_\QQ$ (Lemma \ref{lem:on-the-cone}); and analytically, as a section of $\cH^c_T$ (Lemma \ref{lem:define-a-function}).

\begin{lemma}\label{lem:on-the-cone}
If we expand each $y_i^{\pm} = \exp(\log(y_i^{\pm}))$ appearing in $I^c_T$ as a power series in $\log(y_i^{\pm})$, then the seies $I^c_T$ is a point of the Lagrangian cone of $X\sslash_{\theta_c} T$ valued in the topological $\Lambda^c_{F_K}$-algebra
\[
\Lambda^c_{F_K} \llbracket \log(y_1^-), \ldots, \log(y_c^-), \log(y_{c+1}^+), \ldots, \log(y_k^+) \rrbracket.\]
\end{lemma}
% \begin{remark}\label{rmk:on-the-cone}Since it agrees with the series in CITE CIJ], the series $I^c_T$ is on the Lagrangian cone of $X\sslash_{\theta_c} T$. We view $I^c_T$ as an $R^c$-valued point of the cone, where \[
% R^c = \Lambda^c \llbracket \log(y_1^-), \ldots, \log(y_c^-), \log(y_{c+1}^+), \ldots, \log(y_k^+), x_1, \ldots, x_M \rrbracket,\] by expanding each $y_i^{\pm} = \exp(\log(y_i^{\pm}))$ appearing in $I^c_T$ as a power series in $\log(y_i^{\pm})$.
% \end{remark}
\begin{proof}
This follows from Remark \ref{rmk:agrees} and \cite{CCIT15, ciocan2010moduli, wcgis0}.

\end{proof}

\begin{lemma}\label{lem:define-a-function}
If we set $q_i^{\pm 1}=1$ in the series $I^{c}_T$, we obtain a multivalued section of $\cH^{c}_T$ defined on the open subset $D_c$ of $ U_c$ given by
\begin{equation}\label{eq:define-a-func}
 D_c := \{\; (y_1^-, \ldots, y_c^-, y_{c+1}^+, \ldots, y_k^+) \in U_c \;\; \mid \;\; 0 < |y_i^{\epsilon^c_i}| < 1 \;\}.
\end{equation}
\end{lemma}
\begin{remark}\label{rmk:define-a-section}
The series $I^c_T$ is a multivalued section of $\cH^c_T$ due to the presence of $\log(y_i^{\epsilon^c_i})$ variables. If we remove a ray from 0 to infinity in each $y_i^{\epsilon^c_i}$-plane and choose a branch of each $\log(y_i^{\epsilon^c_i})$, we obtain a genuine section of $\cH^c_T$ on this smaller domain.
\end{remark}

\begin{proof}[Proof of Lemma \ref{lem:define-a-function}]
If we define formal variables $\alpha_{ij} := (H_i-\lambda_j)/z$ and $\beta_{ij}:= (-H_i+\sigma_j)/z$, then we can rewrite the series $I^c_T$ as 
\begin{equation}\label{eq:I-shape}
I^c_T|_{q_i=1} 
= ze^{\square_c/z}  
\prod_{i=1}^c J_i^{-} \prod_{i=c+1}^k J_i^+
\end{equation}
where we set
\[
J_i^+ =  \sum_{d\geq 0} (y_i^+)^d \prod_{j=1}^{n}\left( {\prod_{\ell=1}^{d }(\alpha_{ij}+\ell)^{-1}} \prod_{\ell=-d+1}^0 (\beta_{ij}+\ell) \right),\]
\[
J_i^- = \sum_{d\leq 0} (y_i^-)^{-d} \prod_{j=1}^{n}\left( {\prod_{\ell=d+1}^{0 }(\alpha_{ij}+\ell)} \prod_{\ell=1}^d (\beta_{ij}+\ell)^{-1} \right).
\]
The exponential factor $e^{\square_c/z}$ defines a multivalued section of $\cH^c_T$ on the open set $D_c$, so it is enough to show that each $J_i^{\pm}$ has a radius of convergence equal to 1---ie., the coefficient in $J_i^{\pm}$ of any monomial in the $\alpha_{ij}$ and $\beta_{ij}$ is a series in $y_i^{\pm}$ whose radius of convergence is 1. For this, we use Lemma \ref{lem:flip-flopII} below: if we let $\alpha_{ij}$ and $\beta_{ij}$ be complex numbers of norm less than 1,\footnote{Requiring the norms to be less than 1 ensures that when we substitute these numbers into the formula for $J_i^{\pm}$, we will not be dividing by zero.} one can use the ratio test to show that the resulting series converges when $|y_i^+|<1$ (or $|y_i^-| < 1$).  In the case $i>c$, the relevant limit is
\begin{align*}
\lim_{d \to \infty} &
\prod_{j=1}^n
\frac{\prod_{\ell=-d}^0|\ell+\beta_{ij}|}{\prod_{\ell=1}^{d+1} |\ell+\alpha_{ij}|} 
\frac{\prod_{\ell=1}^{d} |\ell+\alpha_{ij}|}{\prod_{\ell=-d+1}^0 |\ell+\beta_{ij}|}|y_i^+|\\
&= \lim_{d \to \infty} \prod_{j=1}^n\frac{|-d+\beta_{ij}|}{ |d+1 +\alpha_{ij}|} |y_i^+| 
\end{align*}
% \begin{align*}
% \lim_{d \to \infty} &
% \prod_{j=1}^n
% \frac{\prod_{\ell=-d}^0|\ell+\gamma_j-\beta|}{\prod_{\ell=1}^{d+1} |\ell-\alpha_j+\beta|} 
% \frac{\prod_{\ell=1}^{d} |\ell-\alpha_j+\beta|}{\prod_{\ell=-d+1}^0 |\ell+\gamma_j-\beta|}|y_r^+|\\
% &= \lim_{d \to \infty} \prod_{j=1}^n\frac{|-d+\gamma_j-\beta|}{ |d+1 -\alpha_j+\beta|} |y_r^+| 
% \end{align*}
which equals $|y_i^+|$.
\end{proof}

\begin{lemma}\label{lem:flip-flopII}
Let $\bx = x_1, \ldots, x_n$ be complex numbers and suppose that $\phi(\bx) = \sum_{I \in \NN^n} a_I\bx^I$ converges absolutely for some complex coefficients $a_I$. Then for a given index $J \in \NN^{n-1}$ the coefficient $\sum_{m\geq 0} a_{(J, m)}x_n^m$ of $x_1^{J_1}\ldots x_{n-1}^{J_{n-1}}$ in $\phi(\bx)$ converges absolutely.
\end{lemma}
\begin{proof}
We have
\[
\sum_{m} |a_{(J, m)}| |x_n|^m = \frac{\sum_m |a_{(J, m)}| |x_1|^{J_1} \cdots |x_{n-1}|^{J_{n-1}}|x_n|^m}{|x_1|^{J_1} \cdots |x_{n-1}|^{J_{n-1}}}
 \leq \frac{\sum_{J, m}|a_{(J, m)}| |x_1|^{J_1} \cdots |x_{n-1}|^{J_{n-1}}|x_n|^m}{|x_1|^{J_1} \cdots |x_{n-1}|^{J_{n-1}}}.\]
Since $\phi(\bx)$ converges absolutely, the right hand side is a finite sum, so the left hand side is as well.
\end{proof}

\subsection{Step 3: Analytic continuation of small toric $I$-functions}

%If $V \subset \cM_T$ is a locally closed analytic subvariety excluding $p_0$ and $p_k$, we define \textit{a path in $V$ connecting $p_0$ and $p_k$} to be the continuous image of a parameter $t \in (0, 1)$, such that the limit (in the analytic topology) as $t \to 0$ is $p_0$ and the limit as $t \to 1$ is $p_k$. 
Set $I^+_T := I^0_T$ and $I^-_T := I^k_T$, and analogously define $H^+_T, H^-_T$ and their associated sheaves. 
% Finally, if $\fI \subset \RR$ is an interval of length less than $2\pi$ and $y$ is any complex variable, let $\log_{\fI}(y)$ denote the unique value of $\log(y)$ with imaginary part in $i\fI$. We fix the intervals
% \[
% \fI = (i(n-2)\pi , in\pi) \quad \quad \quad \quad -\fI = (-in\pi, i(2-n)\pi).
% \]

\begin{proposition}\label{prop:U-exists}
There is a path $\gamma: [0,1] \to \cM_G$ and an $F_K((z^{-1}))$-module isomorphism $\UU_T: H^+_T \to H^-_T$ satisfying
\begin{itemize}
\item $\gamma(0)$ is in $D_0$ and $\gamma(1)$ is in $D_k$,
\item $\UU_T$ is symplectic, preserves degree, and is $W$-equivariant, and
\item $\UU_T I^+_T|_{q_i=1}$ analytically continues to $I^-_T|_{q_i=1}$ along $\gamma$.
\end{itemize}
\end{proposition}

We will prove Proposition \ref{prop:U-exists} in Section \ref{sec:proof} below: the basic idea is to compose the analytic continuations in \cite{coates2018crepant} and then homotope the resulting path to lie in $\cM_G$. Before giving the proof, we recall the specifics of the analytic continuations constructed in \cite{coates2018crepant} (Section \ref{sec:cij}) and we use a differential equation to construct analytic continuations of the series $I^c_T|_{q_i=1}$ (Section \ref{sec:an-cont}).

\subsubsection{Review of analytic continuations in \cite{coates2018crepant}}\label{sec:cij}
For $c=1, \ldots, k$ let $M_c$ be the (complex analytic) toric subvariety of $\cM_T$ corresponding to the hyperplane wall separating $\theta_{c-1}$ from $\theta_c$ in $\chi(T)_\QQ$. The variety $M_c$ is isomorphic to $\PP^1$ and contains $p_{c-1}$ and $p_c$ as the torus fixed points. (For an example when $k=2$ see Figure \ref{fig:only}.) The variables $y_c^{\pm}$ are natural coordinates on $M_c$, and in fact give a canonical identification of $M_c$ with $\PP^1$. We define a formal analytic lift $\cM_c$ of $M_c$ by assigning the sheaf of functions on $\cM_c$ to be 
\[
\cO_{\cM_c} := \cO_{M_c} \otimes_\CC \CC\llbracket y_1^-,\ldots, y_{c-1}^-, y_{c+1}^+,\ldots, y_{k}^+\rrbracket
\]
The relationship between $\cO_{\cM_c}$ and $\cO_{\cM_T}$ is summarized in the following remark.
\begin{remark}\label{rmk:injective}
Let $\iota_c: \cM_c \to \cM_T$ denote the inclusion. There is a natural restriction morphism $\iota_c^{-1} \cO_{\cM_T} \to \cO_{\cM_c}$ sending an analytic function in $y_i^{\epsilon^c_i}$ to its power series expansion centered at zero in the variables $y_i^{\epsilon^c_i}$ for $i\neq c$. This restriction morphism is injective: if two analytic functions on $(\PP^1)^k$ have the same power series expansion at a point, then they are equal in some neighborhood of that point, hence define the same section of the inverse image sheaf. 
\end{remark}

Observe that $\cM_c$ is the formal $\PP^1$ considered in \cite[(5.9)]{coates2018crepant} that corresponds to crossing the hyperplane wall separating $\theta_{c-1}$ and $\theta_c$. For $\delta=1$ or $0$ the series $I^{c-\delta}_T|_{q_i=1}$ defines a section of
\begin{equation}\label{eq:bundle}
\cO_{\cM_c} \otimes H^{c-\delta}_T[ \log(y_1^-),\ldots,\log(y_{c-1}^-),\log(y_{c+1}^+),\ldots,\log(y_k^+) ] 
\end{equation}
in the domain $|y_c^+| < 1$ (if $\delta=1$) or $|y_c^-| < 1$ (if $\delta=0$) by the proof of Lemma \ref{lem:define-a-function}. By this, we mean that the coefficient of any monomial in the $\log(y_i^{\pm})$ (for $i\neq c$) is a section of $\cO_{\cM_c}\otimes H^{c-\delta}_T$.  
% \begin{equation}\label{eq:Ic-1-minus-log}
% e^{-\square_{c-1}/z} I^{c-1}_T|_{q_i=1} 
% % e^{-\square_{c-1}/z + H_c\log(y_c^+))/z}I^{c-1}_T|_{q_i=1} 
% \end{equation}
%  defines a section of $\cO_{\cM_c} \otimes H^{c-1}_T $ in the domain $|y_c^+| < 1$ by the proof of Lemma \ref{lem:define-a-function}.  
%  Likewise
%  \begin{equation}\label{eq:Ic-minus-log}
%  e^{-\square_c/z}I^c_T|_{q_i=1}
% % e^{-\square_c/z - H_c\log(y_c^-))/z}I^c_T|_{q_i=1}
% \end{equation}
% defines a section of $\cO_{\cM_c} \otimes H^{c}_T$ in the domain $0 < |y_c^-| < 1$.
 The key result we use from \cite{coates2018crepant} is that the series $I_T^{c-\delta}$, $\delta=0,1$ are related by analytic continuation (meaning that the coefficients of the
monomials in $\log(y_i^{\pm})$ are related by analytic continuation).
\begin{theorem}[{\cite[Thm 6.1]{coates2018crepant}}]\label{thm:cij}Fix a copy of $\PP^1$ with coordinates $[y^+:y^-]$.
There is a path $\gamma:[0,1] \to \PP^1$ and for each $c=1, \ldots, k$ there is an $F_K((z^{-1}))$-module isomorphism $\UU_c:H^{c-1}_T \to H^c_T$ satisfying
\begin{itemize}
\item $\gamma(0) = [\epsilon:1]$ and $\gamma(1) = [1:\epsilon]$ for some small positive real $\epsilon$ satisfying $0< \epsilon < 1$,
\item $\UU_T$ is symplectic and preserves degree, and
\item After identifying $M_c$ with our fixed $\PP^1$ via the coordinates $y^{\pm}_c$, the series $\UU_cI^{c-1}_T|_{q_i=1}$ analytically continues to $I^c_T|_{q_i=1}$ along $\gamma$ as sections of \eqref{eq:bundle}.
\end{itemize}
\end{theorem}

\begin{proof}
By Remark \ref{rmk:agrees} our series are compatible with the ones used in \cite{coates2018crepant}. From \cite[Thm 6.1]{coates2018crepant}
we have a degree-preserving symplectic isomorphism $\UU_c: H^{c-1}_T \to H^c_T$ and a path $\gamma_c$ in the $\log(y_c^+)$-plane (pictured in \cite[Fig 1]{coates2018crepant}) such that $I^{c-1}_T|_{q_i=1}$ agrees with $\UU_c I^c_T|_{q_i=1}$ after analytic continuation along $\gamma_c$. The only thing we need to check is that $\gamma_c$ can be chosen independent of $c$, but this follows from the definition in \cite[Fig 1]{coates2018crepant}, noting that the constant $w$ there is equal to $n-1$ and $|\mathfrak{c}|$ is equal to 1 for every $c$.
\end{proof}

\subsubsection{Analytic continuations of $I^c_T$}\label{sec:an-cont}

Recall from the proof of Lemma \ref{lem:define-a-function} the expression \eqref{eq:I-shape} for $I^c_T$ in terms of series $J_i^{\pm}$. 
We now view $J_{ i}^{\pm}$ as a function $J_{ i}^{\pm}(y_i^+, \balpha_i, \bbeta_i)$ of $1+2n$ complex variables, where $\balpha = (\alpha_{i1}, \ldots, \alpha_{in})$ and $\bbeta = (\beta_{i1}, \ldots, \beta_{in})$.
Our strategy is to analytically continue $I^c_T|_{q_i=1}$ by analytically continuing the coefficients in $J_{i}^{\pm}$ of monomials in $\alpha_{ij}$ and $\beta_{ij}$. 

\begin{lemma}\label{lem:J}
The function $J^{\pm}_{ i}(y_i^\pm, \balpha, \bbeta)$ is annihilated by a differential operator of the form
\begin{equation}\label{eq:ode}
(d/dy_i^\pm)^{n+1} + a_1(d/dy_i^\pm)^{n} + \ldots + a_{n}(d/dy_i^\pm) + a_{n+1}
\end{equation}
where each $a_i = a_i(y_i^\pm, \balpha_i, \bbeta_i)$ is analytic on $(\CC \setminus \{0, (-1)^n\}) \times \CC^{2n}$.
\end{lemma}
\begin{proof}
We prove the lemma for $J_i^+$.
One can check directly that $J_{i}^+$ satisfies the equation
\begin{equation}\label{eq:gkz}
K+\prod_{j=1}^{n}(\alpha_{ij} J_i^+ + y_i^+\frac{dJ_{i}^+}{dy_i^+}  ) = y_i^+ \prod_{j=1}^n (\beta_{ij} J_i^+ - y_i^+ \frac{dJ_{ i}^+}{dy_i^+} ) \quad \quad \quad \text{where} \quad \quad \quad K = (-1)\prod_{j=1}^{n}
(\alpha_{ij}).
\end{equation}
Indeed, for any complex number $\gamma$ we have $(\gamma+(y_i^+) d/dy_i^+ ) (y_i^+)^d = (\gamma+d) (y_i^+)^d$, so a series $\sum_{d\geq 0} a_d (y_i^+)^d$ satisfies \eqref{eq:gkz} if and only if
\[
K + \sum_{d \geq 0} \prod_{j=1}^n(\alpha_{ij}+d) a_d(y_i^+)^d = y_i^+\sum_{d \geq 0} \prod_{j=1}^n (\beta_{ij}-d)a_d (y_i^+)^d;
\]
i.e., if and only if 
\[
a_0 = -K \left(\prod_{j=1}^n \alpha_{ij} \right)^{-1} \quad \quad \text{and} \quad \quad a_d =  \frac{\prod_{j=1}^n (\beta_{ij} - (d-1))}{\prod_{j=1}^n(\alpha_{ij}+d)} a_{d-1}\quad \quad \text{for}\;d \geq 1.
\]
One checks that the coefficients of $J_i^+$ satisfy this initial condition and recursion. It is straightforward to see that \eqref{eq:gkz} can also be written as
\[
(y_i^+)^n(y_i^+-(-1)^n)(d/dy_i^+)^n J_{ i}^+ + L\cdot J_{ i}^+ + K = 0
\]
where $L$ is a differential operator of order $n-1$ analytic in $y_i^+, \balpha_i$, and $\bbeta_i$. Applying $d/dy_i^+$ to the above equation, we see that $J_{ i}^+$ is annihilated by an operator of the form
\[
(y_i^+)^n(y_i^+-(-1)^n)(d/dy_i^+)^{n+1} + L'
\]
where $L$ has order $n$. Dividing through by $(y_i^+)^n(y_i^+-(-1)^n)$ we obtain the result.
\end{proof}

Let $0 \in M_c$ be the point where $y_c^+=0$ and let $\infty \in M_c$ be the point where $y_c^-=0$. Observe that $\cM_T = M_1 \times \ldots \times M_k$ and that we may write $
p_c = (\infty,  \ldots, \infty, 0, \ldots, 0)$ where the first $c$ coordinates are equal to $\infty$ and the last $k-c$ coordinates are equal to $0$. 
%\Rachel{added this to clarify 0 and $\infty$ notation. Does it look ok?}\Wendelin{Yep}\Qaasim{Yep}

\begin{corollary}\label{cor:an-cont-J}
The coefficient of any monomial in $\alpha_{ij}$ and $\beta_{ij}$ in $J^\pm_i$ can be analytically continued along any path $\gamma$ in 
\[
V_i := \cM_i \setminus\{0, (-1)^n, \infty\}.
\]
with $\gamma(0)$ in the domain of convergence of $J^{\pm}_i$. Moreover, if $\gamma$ and $\gamma'$ are two such paths with $\gamma(1)=\gamma'(1)$ such that $\gamma$ is homotopic to $\gamma'$ via a homotopy lying entirely in $V_i$, then the analytic continuations along $\gamma$ and $\gamma'$ agree.
\end{corollary}
%\Qaasim{Is this corollary trivial if $i \neq c$? and what does the statement mean for $J_{c}^-$?}\Rachel{I changed $c$ to $i$. Does this section make sense now? For the second question, I'm not sure what the problem is.}
\begin{proof}
By restricting \eqref{eq:ode} to the subspace $\balpha_i = \bbeta_i = 0$, we see that the corollary holds for the constant term $J_i^{\pm}(y_i^{\pm},0,0)$. It follows from e.g. \cite[Chapter 1 Thm 8.3]{CodLev} that the analytic continuation of $J_i^{\pm}(y_i^{\pm},0,0)$ is actually a restriction of some $\widetilde{J_{ i}^\pm}(y_i^\pm, \balpha_i, \bbeta_i)$ that is analytic for all $(\balpha_i, \bbeta_i)$ in some neighborhood of 0, hence all coefficients of $\alpha_{ij}$ and $\beta_{ij}$ monomials admit analytic continuations.
\end{proof}
%\Wendelin{Rereading this section as a whole, I was just wondering how much of 5.3.6 (1) follows directly from 5.3.3 (analytic continuation in CIJ)} \Rachel{I think the key in 5.3.6(1) is ANY PATH in $V_1 \times \ldots V_k$. If we just quote CIJ I think we just have analytic continuation on one specific path, or maybe on a smaller disk.}
\begin{corollary}\label{cor:an-cont}
Let $c \in \{1, \ldots, k\}$. We can analytically continue $I^{c-1}_T|_{q_i=1}$ in the following two situations:
\begin{enumerate}
\item As a multivalued section of $\cH^{c-1}_T$, it can be analytically continued along any path contained in $V_1 \times \ldots \times V_k \subset \cM_T$ with $\gamma(0) \in D_{c-1}$.
\item As a multivalued section of 
\begin{equation}\label{eq:funny sheaf}
\cO_{\cM_T} \otimes H^{c-1}_T[\log(y^{-}_1), \ldots, \log(y^-_{c-1}), \log(y^+_{c+1}), \ldots, \log(y_k^+)]
\end{equation}
it can be analytically continued along any path contained in $\infty \times\ldots \times \infty \times V_{c} \times 0 \ldots \times 0 \subset \cM_{c} \subset \cM_T$ with $\gamma(0) \in D_{c-1}$. 
%We just consider the $\log(y_i^{\pm})$ as formal variables here.
\end{enumerate}
In either case, the analytic function obtained near $\gamma(1)$ depends only on the homotopy class of the path, provided the homotopy is contained in the specified domain.
\end{corollary}

% \begin{corollary}\label{cor:an-cont}
% We can analytically continue $I^{c}_T|_{q_i=1}$ along the following two kinds of paths:
% \begin{enumerate}
% \item Any path contained in $V_1 \times \ldots \times V_k \subset \cM_T$ with $\gamma(0) \in D_c$.
% \item Any path contained in $\infty \times\ldots \times \infty \times V_c \times 0 \ldots \times 0 \subset \cM_c \subset \cM_T$ with $\gamma(0) \in D_c$. \Rachel{probably something is missing in here: perhaps analytically continue coefficients of log variables and argue that expansion in terms of log variables is unique.}
% \end{enumerate}
% In either case, the analytic function obtained near $\gamma(1)$ depends only on the homotopy class of the path, provided the homotopy is contained in the specified domain.
% \end{corollary}
\begin{proof}
These analytic continuations are constructed by taking products of the analytic continuations in Corollary \ref{cor:an-cont-J}. In the second part, the key point is that even though the 0's and $\infty$'s appearing in coordinates $i$ for $i\neq c$ are outside the domain $V_i$ of analytic continuation, the path is constant in these coordinates and so we may use the constant ``analytic continuation'' of coefficients of $J^{\pm}_i$.
\end{proof}

\subsubsection{Proof of Proposition \ref{prop:U-exists}}\label{sec:proof}

Let $\gamma': [0,1] \to \PP^1$ be the path in Theorem \ref{thm:cij} and let $\UU_1, \ldots, \UU_k$ be as in that theorem. Define
\[
\gamma(t) = (\gamma'(t), \ldots, \gamma'(t)) \quad \quad \quad \quad \UU_T = \UU_k \UU_{k-1}\ldots \UU_1
\]
It follows from Corollary \ref{cor:an-cont} that $I^+_T|_{q_i=1}$ admits an analytic continuation along $\gamma$.

We now show that $\UU_T I^+_T|_{q_i=1}$ analytically continues to $I^-_T|_{q_i=1}$ along $\gamma$.
To do so, we first construct an analytic continuation along an alternative path. By Corollary \ref{cor:an-cont} we have that $I^{c-1}_T|_{q_i=1}$ can be analytically continued as a multivalued section of \eqref{eq:funny sheaf} along the path
\[
(\infty, \ldots, \infty, \gamma'(t), 0, \ldots, 0) \quad \quad \subset \cM_c.
\]
But this is precisely the path appearing in \cite{coates2018crepant}, and it follows from Theorem \ref{thm:cij} and Remark \ref{rmk:injective} together with uniqueness of analytic continuations that $\UU_cI^{c-1}_T|_{q_i=1}$ agrees with $I^{c}_T|_{q_i=1}$ after analytic continuation along this path, at least as sections of \eqref{eq:funny sheaf}. In fact, by shrinking the constant $\epsilon$  appearing in Theorem \ref{thm:cij}, we can assume that $\UU_cI^{c-1}_T|_{q_i=1}$ analytically continues to $I^{c}_T|_{q_i=1}$ along the path
\[
\gamma'_c(t):=([1:\epsilon],\ldots [1:\epsilon], \gamma'(t), [\epsilon:1], \ldots, [\epsilon:1]) \quad \quad \subset \cM_T.
\]
In a neighborhood of this new path, we have a morphism sending sections of \eqref{eq:funny sheaf} to multivalued sections of $\cM_T$, induced by viewing the formal variables $\log(y_i^{\pm})$ for $i \neq c$ as multivalued analytic functions.
%(\Wendelin{the variables are constant only on the path, not a neighbourhood}), but multivalued, analytic functions. \Rachel{try this instead: ``a morphism sending sections of \eqref{eq:funny sheaf} to multi-sections of $\cM_T$, induced by viewing the formal variables $\log(y_i^{\pm})$ for $i\neq c$ as multivalued analytic functions. The key here is that since $\epsilon \neq 0$, in a neighborhood of $[1:\epsilon]$ the possible values of $y_i^-$ are a set of single-valued analytic functions.'' should we talk about multi-sections rather than multi-valued sections? a multi-section would be a set of sections}\Wendelin{I think multi-valued is fine.} 
It follows that $\UU_cI^{c-1}_T|_{q_i=1}$ analytically continues to $I^{c}_T|_{q_i=1}$ along $\gamma'_c(t)$ as multivalued sections of $\cM_T$, and hence that $\UU_T  I^{+}_T|_{q_i=1}$ analytically continues to $I^{-}_T|_{q_i=1}$ along the concatenation of paths $\gamma'_k \star \ldots \star \gamma'_1$. Since this concatenation lies in $V_1 \times \ldots \times V_c$, all we have to do is give a homotopy from $\gamma$ to $\gamma'_k \star \ldots \star \gamma'_1$ that lies entirely in $V_1 \times \ldots \times V_c$. If we let $p_i^\epsilon \in \cM_T$ denote the point whose first $i$ coordinates are equal to $[1:\epsilon]$ and whose remaining coordinates are equal to $[\epsilon:1]$, then one such homotopy is given by $h_s(t) = (h^1_s(t), \ldots h^k_s(t))$, where
\begin{equation}\label{eq:homotopy}
h^i_s(t) = \left\{\begin{array}{ll}
p_{i-1}^\epsilon & 0 \leq t \leq sk^{-1}(i-1) \\
\gamma \left(\frac{kt-s(i-1)}{k - s(k-1)} \right) & sk^{-1}(i-1) \leq t \leq 1-sk^{-1}(k-i) \\
p_i^\epsilon & 1-sk^{-1}(k-i) \leq t \leq 1
\end{array}\right.
\end{equation}
and we note that $h_0(t) = \gamma(t)$ and $h_1(t) = \gamma'_k \star \ldots \star \gamma'_1.$ %\Rachel{slightly more details in details.tex} 
This completes the proof that $\UU_T I^+_T|_{q_i=1}$ analytically continues to $I^-_T|_{q_i=1}$ along $\gamma.$

To complete the proof of Proposition \ref{prop:U-exists}, we have to show that $\UU_T$ is $W$-equivariant. To make the argument completely transparent, we articulate our understanding of analytic continuation in the following remark.
\begin{remark}[On the meaning of analytic continuation] \label{rmk:meaning}
We have shown that $\UU_T I^+_T|_{q_i=1}$ analytically continues to $I^-_T|_{q_i=1}$ along $\gamma$; we express this relationship symbolically as
\[
\UU_T I^+_T|_{q_i=1} \simeq_\gamma I^-_T|_{q_i=1}.
\]
The relationship $\simeq_\gamma$ means that we have an increasing sequence of real numbers 
\[
0=b_0 < b_1 < \ldots < b_m=1,
\]
a neighborhood $B_\ell$ of $\gamma(b_\ell)$ for each $\ell=0, \ldots, m$, and a section $\widetilde{I^+_{T, \ell}}\in \cH^+_T(B_\ell)$ such that $\gamma$ maps $[b_\ell, b_{\ell+1}]$ into $B_\ell \cup B_{\ell+1}$ for $\ell=0, \ldots, m-1$, we have that $\widetilde{I^+_{T, \ell}} = \widetilde{I^+_{T, \ell+1}}$ on $B_\ell \cap B_{\ell+1}$, and we have $\widetilde{I^+_{T, 0}} = I^+_T|_{q_i=1}$ and
\begin{equation}\label{eq:with-y}
\UU_T \widetilde{I^+_{T, m}}|_{y_i^+ = (y_i^-)^{-1}} = I^-_T|_{q_i=1}
\end{equation}
%\Wendelin{LHS of eq 24 should not need the $q_i=1$, right?}
where the domains of these sections overlap.
Since $I^+_T$ uses the coordinate $y^+_i$, we think of each $\widetilde{I^+_{T, \ell}}$ as written in the coordinate $y^+_i$ and hence we articulate the change of variables in \eqref{eq:with-y}.
\end{remark}
\begin{lemma}\label{lem:unique}
The transformation $\UU_T: H^+_T \to H^-_T$ is the unique $F_K((z^{-1}))$-linear map such that $\UU_T I^+_T|_{q_i=1} \simeq_\gamma I^-_T|_{q_i=1}$.
\end{lemma}
\begin{proof}
We remark that this argument is used in \cite{coates2018crepant}; see, for example, \cite[(5.17)]{coates2018crepant}.

Consider the set of monomials
\[
S = \{\mu(u_1, \ldots, u_k) = u_1^{b_1} \cdots u_k^{b_k} \mid 0 \leq b_i \leq n-1 \}.
\]
Observe that the set $\{\mu(-H_1, \ldots, -H_k) \mid \mu \in S\}$ of Chow classes is a basis for $A^*_K(X\sslash_{\theta_-} T)_\QQ$, and hence for the $F_K((z^{-1}))$-vector space $H^-_T$: this set is the usual basis for the equivariant Chow ring of $(\PP^{n-1})^k$, and $X\sslash_{\theta_-} T$ is a vector bundle over this space. On the other hand, if we define differential operators
\[
\partial_{\log(y^-)}^\mu := \mu(z\partial_{\log(y_1^-)}, \ldots, z\partial_{\log(y_k^-)})
\]
then we claim that if $p$ is any point in the domain of convergence of $I^-_T$, the derivatives $\partial_{\log(y^-)}^\mu I^-_T$ evaluated at $p$ are linearly independent elements of $H^-_T$. Since differentiation commutes with the linear map $\UU_T$, this claim implies the lemma. For the claim, a direct computation shows 
\[
\partial_{\log(y^-)}^\mu I^-_T = e^{\square_k/z}(\mu(-H_1, \ldots, -H_k) + O(y_i^-)),
\]
so that there is an invertible matrix with entries in the ring $F_K((z^{-1}))[\log(y_i^-)]\llbracket y_i^-\rrbracket$ whose columns are the series $\partial_{\log(y^-)}^\mu I^-_T$ written in the basis defined by $S$. 
%\Qaasim{why does this follow from the above computation?} \Rachel{Does this make sense now? the matrix is invertible because the inverse of $e^{\square_k/z}$ is $e^{-\square_k/z}$ and the second factor has the form $I + O(y)$ where $I$ is the identity matrix and $O(y)$ is a matrix whose entries are in $O(y)$ (no log variables, potentially z variables). This $O(y)$ is topologically nilpotent where the topology is just about the $y$ variables. The inverse of $I + O(y)$ can be ``computed'' using the geometric series: $1/(I+O(y)) = 1 - O(y) + O(y)^2 - \ldots$ and it is important that $O(y)$ is nilpotent/doesn't contain constant terms so that the sum converges.}\Qaasim{Thanks for this, I agree this shows that the matrix is invertible. I'm still trying to understand why this implies the result. I would like to say that this provides a matrix with linearly independent columns whose entries are in $H^-_T$, but I'm still confused.} \Rachel{yes, plug in the point $p$ now mentioned earlier in the proof}\Qaasim{Ahh yes sorry you have said this above, happy now.}
\end{proof}

We apply the lemma to show that $\UU_T$ is Weyl-equivariant. 
Recall that the Weyl group $S_k$ acts on $A^*_K(X\sslash_{\theta_\pm}T)_\QQ$. For $w$ a permutation we set $w \cdot H_i = H_{w(i)}$, and we let $w I^{\pm}_T$ denote the series obtained from letting $w$ act on the Chow classes in $I^{\pm}_T$. A direct computation shows
\begin{equation}\label{eq:b1}
w \,I^{\pm}_T\,|_{y_i^{\pm} = y_{w(i)}^{\pm},\; q_i=1} = I^{\pm}_T|_{q_i=1}.
\end{equation}
It follows by induction and uniqueness of analytic continuations that
\begin{equation}\label{eq:b2}
w \,\widetilde{I^+_{T, \ell}}\,|_{y_i^+ = y_{w(i)}^+} = \widetilde{I^+_{T,\ell}}.
\end{equation}
for each $\ell$, where $\widetilde{I^+_{T, \ell}}$ was defined in Remark \ref{rmk:meaning}. Combining \eqref{eq:with-y} with \eqref{eq:b1} and \eqref{eq:b2} yields
\[
(w^{-1}\UU_T w)\; \widetilde{I^+_{T,m}}\,|_{y_i^+ = (y_{w(i)}^-)^{-1}} = I^-_T|_{y_i^{-} = y_{w(i)}^{-},\; q_i=1}
\]
and applying the change of variables $y^-_i \mapsto y^-_{w^{-1}(i)}$ yields
\[
w^{-1}\UU_T w \widetilde{I^+_{T,m}}|_{y_i^+ = (y_{i}^-)^{-1}} = I^-_T|_{ q_i=1}.
\]
It follows from Lemma \ref{lem:unique} that $w^{-1}\UU_T w = \UU_T$.

\subsection{Step 4: Analytic continuation of big toric $I$-functions}

Define a set of monomials
\[
\{ \mu(u_1, \ldots, u_k) = u_1^{b_1} \cdots u_k^{b_k} \mid 0 \leq b_i \leq n-k \;\text{and}\; \sum_{i=1}^k b_i \neq 1\}.
\]
Let $M$ denote the cardinality of this set and enumerate its elements $\mu_1, \ldots, \mu_M$. For any $\mu_j$ define
\begin{equation}\label{eq:plug-mu}
\mu_j(\bH + \bd z) = \mu_j(H_1 + d_1z, \ldots, H_k+d_kz).
\end{equation}
Finally, let $x_1, \ldots, x_M$ be formal variables, set $\square_+ := \square_0$ and $\square_- = \square_k$, and define
\[
\II^{\pm}_T := ze^{\square_{\pm}/z}\sum_{\pm d_i \geq 0} \prod_{i=1}^k q_i^{d_i} (y_i^{\pm})^{\pm d_i} e^{\sum_{j=1}^M x_j \mu_j(\bH + \bd z)/z}  I^{\pm}_{T, \bd}
\]
Observe that $\II^{\pm}_T |_{x_i=0} = I^{\pm}_T$.

\begin{remark}
It follows from either \cite{bigI} or \cite{giv15} that $\II^{\pm}_T$ is on the Lagrangian cone of $X\sslash_{\theta_{\pm}}T$, but we will not use this.

\end{remark}

The key insight is that the series $\II^{\pm}_T$ are obtained by differentiating $I^{\pm}_T$ and hence they inherit the analytic properties of $I^{\pm}_T$. Set $D_+ =D_0 $ and $D_- = D_k$ where $D_c$ is the domain defined in \eqref{eq:define-a-func}.

\begin{lemma}\label{lem:big-I}
The series $\II^{\pm}_T|_{q_i=1}$ is a multivalued section of \[\cO_{\cM_T}\llbracket x_1, \ldots, x_M \rrbracket \otimes_\CC H^{\pm}_T\]
on the domain $D_{\pm}$.
Moreover, Proposition \ref{prop:U-exists} holds with $\II^+_T$ and $\II^-_T$ in place of $I^+_T$ and $I^-_T$, respectively (but with the same $\gamma$ and $\UU_T$ as before).
\end{lemma}
\begin{proof}
The series $\II^+_T|_{q_i=1}$ is equal to
\[\sum_{d_i \geq 0} \sum_{\ba \in \NN^M} \frac{(x_1\mu_1(\bH+\bd z)/z)^{a_1} \ldots (x_M\mu_M(\bH+\bd z)/z)^{a_M}}{a_1!\ldots a_M!}e^{\sum_{i=1}^k \log(y_i^{+})( H_i+ d_iz)/z} I^{\pm}_{T,\bd}.\]
%&=\sum_{d_i \geq 0} \sum_{\ba \in \NN^M} \frac{\bx^{\ba}}{\ba ! z^{|\ba|}} ((H_1+D_1z))^{(P\ba)_1} \ldots ((H_k+D_kz)/z)^{(P\ba)_k}e^{\sum_{i=1}^k \log(y_i^{+})( H_i+ d_iz)/z} I^{\pm}_{T,\bd}\\
If we define
\[
\frac{\bx^{\ba}}{\ba ! z^{|\ba|}} = \frac{x_1^{a_1}\cdots x_M^{a_M}}{a_1!\cdots a_M! z^{a_1+\ldots+a_M}} \quad \quad \text{and} \quad \quad \partial^{\mu^{\ba}}_{\log(y^\pm)} = \prod_{i=1}^M \mu_i(z\partial_{\log(y^\pm_1)}, \ldots, z\partial_{\log(y^\pm_k)})^{a_i},
\]
then we can write
\begin{equation}
\II^+_T|_{q_i=1}= \sum_{\ba \in \NN^M} \frac{\bx^{\ba}}{\ba ! z^{|\ba|}} \partial^{\mu^{\ba}}_{\log(y^+)} {I^+_T}|_{q_i=1}.\label{eq:big-computation}
\end{equation}
%\Qaasim{$I^+_T$ is analytic as a function of the $y_i$, the differentation is with respect to $\log y_i$, that's fine right? I'm thinking about the chain rule but confused about introducing $1/y_i$ ahh no it's fine I think I got the chain rule the wrong way round, it would introduce a $y_i$}It follows that $\II^+_T|_{q_i=1}$ is analytic, and we may analytically continue it along $\gamma$ using the functions
\[
\widetilde{\II^+_{T,\ell}} := \sum_{\ba \in \NN^M} \frac{\bx^{\ba}}{\ba ! z^{|\ba|}} \partial^{\mu^\ba}_{\log(y^+)} \widetilde{I^+_{T,\ell}} \quad \quad \quad \quad \ell=1, \ldots, m
\]
where $\widetilde{I^+_{T,\ell}}$ was defined in Remark \ref{rmk:meaning}. To finish the proof we need to show that
$\UU_T \widetilde{\II^+_{T,m}}|_{y_i^+ = (y_i^-)^{-1}}$ agrees with  $\II^-_T|_{q_i=1}$ on a common domain. If $\mu^{\ba}(-1) = \prod_{i=1}^M \mu_i(-1, \ldots, -1)^{a_i}$, then we have
\begin{align*}
\UU_T \widetilde{\II^+_{T,m}}|_{y_i^+ = (y_i^-)^{-1}} &= \sum_{\ba \in \NN^M} \frac{\bx^{\ba}}{\ba ! z^{|\ba|}} \UU_T\left(\partial^{\mu^\ba}_{\log(y^+)} \widetilde{I^+_{T,m}}\right)|_{y_i^+ = (y_i^-)^{-1}} \\
%&= \sum_{\ba \in \NN^M} \frac{\bx^{\ba}}{\ba ! z^{|\ba|}} \UU_T(-1)^{\sum_i (Q\ba)_i} %\partial^{Q\ba}_{\log(y^-)}\left( \widetilde{I^+_T}|_{y_i^+ = (y_i^-)^{-1}} \right)\\
&= \sum_{\ba \in \NN^M} \frac{\bx^{\ba}}{\ba ! z^{|\ba|}} \mu^{\ba}(-1) \partial^{\mu^\ba}_{\log(y^-)}\left( \UU_T \widetilde{I^+_{T,m}}|_{y_i^+ = (y_i^-)^{-1}} \right)\\
\end{align*}
where we have used the definition of $\widetilde{I^+_{T, m}}$ followed by the chain rule and linearity of differentiation.
Since $\UU_T \widetilde{I^+_{T,m}}|_{y_i^+ = (y_i^-)^{-1}}$ restricts to $I^-_T|_{q_i=1}$ by Proposition \ref{prop:U-exists}, we see that after restriction we have
\[
\UU_T \widetilde{\II^+_{T,m}}|_{y_i^+ = (y_i^-)^{-1}} = \sum_{\ba \in \NN^M} \frac{\bx^{\ba}}{\ba ! z^{|\ba|}} \mu^\ba(-1) \partial^{\mu^\ba}_{\log(y^-)}\left( {I^-_T}|_{q_i=1}\right).
\]
The right hand side is equal to $\II^-_T|_{q_i=1}$ by a computation analogous to the one that produced \eqref{eq:big-computation}: applying $z\partial_{\log(y_i^-)}$ to $I^-_T|_{q_i=1}$ multiplies it by $-(H_i+D_iz)$, and the negative signs cancel with the sign that came out of the chain rule.

\end{proof}

\subsection{Step 5: Analytic continuation of nonabelian $I$-functions}
The Weyl group $W=S_k$ acts on the set $\{\mu_j\}_{j=1}^M$ via its action on the $H_i$, hence it also acts on the set of indices $\{1, \ldots, M\}$. Let $N$ denote the number of orbits, and let
\[
J: \{1, \ldots, M\} \to \{1, \ldots, N\}
\]
be the function that sends an index $i$ to the index of its Weyl-orbit.
Define a set of polynomials $P_1, \ldots, P_N$ via
\[
P_j(u_1, \ldots, u_k) := \sum_{J(i)=j} \mu_i(u_1, \ldots, u_k).
\]
Define $P_j(\bH+\bd z)$ in analogy with \eqref{eq:plug-mu}. Introduce the Novikov ring $\Lambda^\pm_{F_K} := F_K\llbracket \bq^{\pm} \rrbracket$ and let $\bx_1, \ldots, \bx_N$ be formal variables. Recall the basis of characters $\{e_i\}_{i=1}^k$ from Section \ref{sec:step1}. If we set $pr^\pm_i := \pm e_i$ then $\{pr^{\pm}_i\}_{i=1}^k$ are bases of $\chi(T)$, and the associated series \eqref{eq:nonab1} agree with our series $\II^{\pm}_T|_{x_i=\bx_{J(i)}}$. If we let $\Delta, \partial_{\Delta}$, and $\zeta$ be as in Section \ref{sec:abelianization} and write $\zeta = \sum_{i=1}^k a^{\pm}_i pr_i^{\pm},$ then by Theorem \ref{thm:nonabelianI} there are unique series $\II^{\pm}_G$ satisfying 
\begin{equation}\label{eq:computation}
\Delta g^*\II^{\pm}_G =  \partial_\Delta j^*\II^{\pm}_T|_{y^{\pm}_i=\by^{\pm}, \;x_i = \bx_{J(i)},\; q_i= (-1)^{a^\pm_i}\bq}.
\end{equation}

\begin{remark}
The objects $\Delta, \partial_{\Delta},$ and $\zeta$ used in the definition of $\II^{\pm}_G$ depend on a choice of positive roots $\Phi_+$. It is convenient to use the same set of positive roots $\Phi_+ \subset \chi(T)$ to define both $\II^+_G$ and $\II^-_G$ below, even though by Remark \ref{rmk:choice of roots} this choice does not affect the actual definition of the series. 

The bases $\{pr^{\pm}_i\}$ determine how $\partial_\Delta$ looks in coordinates. For example, if we choose the standard set of positive roots $e_i-e_j$ for $1\leq i < j \leq k$, then 
\begin{equation}\label{eq:par-delta}
\partial_\Delta = \prod_{1\leq i < j \leq k} (\partial_{\log(y_i^+)}-\partial_{\log(y_j^+)}) = \prod_{1\leq i < j \leq k} (-\partial_{\log(y_i^-)}+\partial_{\log(y_j^-)}).
\end{equation}

Moreover, using the same set of positive roots, one can check directly that
\begin{equation}\label{eq:simplify-sign}
(-1)^{a^\pm_i} = (-1)^{k-1}.
\end{equation}
\end{remark}

The following lemma gives us two ways to interpret $\II^{\pm}_G$.

\begin{lemma}\label{lem:interpret-IG}
The series $\II^{\pm}_G$ is a big point of the Lagrangian cone of $X\sslash_{\theta_\pm} G$ valued in the topological $\Lambda^{\pm}_{F_K}$-algebra
\[
\Lambda^{\pm}_{F_K}\llbracket \bx_1, \ldots, \bx_N, \log(\by^{\pm}) \rrbracket.
\]
Moreover, the series $\II^{\pm}_G|_{\bq=(-1)^{k-1}}$ is a multivalued section of
\begin{equation}\label{eq:interpret-IG}
\cO_{\cM_G}\llbracket \bx_1, \ldots, \bx_N \rrbracket \otimes_\CC H^{\pm}_G
\end{equation}
on the open subset of $U^{\pm} = \Spec(\CC[\by^{\pm}])$ where $0 < |\by^{\pm}| < 1$.
\end{lemma}
\begin{proof}
That $\II^{\pm}_G$ is on the cone follows from Theorem \ref{thm:nonabelianI}, since the morphism $A_K^*(X\sslash_{\theta_{\pm}}G)_\QQ \to H^*_K(X\sslash_{\theta_{\pm}} G; \QQ)$ is an isomorphism in this case (it is an isomorphism for $Gr(k, n)$ and both Chow and cohomology are homotopy invariant). The proof of Theorem \ref{thm:nonabelianI} together with \eqref{eq:simplify-sign} shows
\begin{equation}\label{eq:def-IG}
\Delta g^*\II^{\pm}_G = z\partial_\Delta \sum_{\pm d \geq 0} (-1)^{d(k-1)}\bq^d(\by^{\pm})^{\pm d} \sum_{\substack{\pm d_i \geq 0\\ \sum_{i=1}^k d_i=d}} e^{\log(\by^{\pm})\sum_{i=1}^k (\pm H_i\pm d_iz)/z + \sum_{j=1}^N \bx_j P_j(\bH + \bd z)/z } j^*I^{\pm}_{T, \bd}.
\end{equation}
From here straightforward computation shows (keeping in mind the isomorphism \eqref{eq:chow-diagram}) %\Rachel{not sure what you want here}\Qaasim{Should $t$ be $\mathbf{x}$ and $\log(\mathbf{y})$?} \Rachel{yes i think so, also bold $\bq$ right?}
\[
\frac{\partial \II^{\pm}_G}{\partial \bx_j} = P_j(\bH) + O(\bq, \bx, \log(\by^{\pm})) \quad \quad \quad \quad \frac{\partial \II^{\pm}_G}{\partial \log(\by^{\pm})} = \pm \sum_{i=1}^k H_i + O(\bq, \bx, \log(\by^{\pm})),
\]
so to show that $\II^{\pm}_G$ is big, according to Definition \ref{def:big} we only have to check that the set
\[
\{\pm\sum_{i=1}^k H_i, P_1(\bH), \ldots, P_N(\bH)\}
\]
is a basis for $A^*_K(X\sslash_{\theta_\pm}G)_\QQ$. But this set is a well-known basis for the equivariant Chow ring of the Grassmannian $Gr(k, n)$ (at least up to the sign on the first basis element), so by homotopy invariance of equivariant Chow groups we obtain the result.

Analyticity of $\II^{\pm}_G|_{\bq=(-1)^{k-1}}$ follows from the expression \eqref{eq:computation} for $\II^{\pm}_G$ as a derivative of $\II^{\pm}_T$ and Lemma \ref{lem:big-I}, which asserts that $\II^{\pm}_T|_{q_i=1}$ is analytic.
\end{proof}

% \Rachel{remove this remark? and the analogous abelian one earlier?}
% \begin{remark}
% The series $\II^{\pm}_G|_{\bq=1}$ is characterized as a multivalued section of \eqref{eq:interpret-IG} in Lemma \ref{lem:interpret-IG} because of the presence of $\log(\by^{\pm})$ in its definition. On any simply connected subset of $\{\by^{\pm}|0 < |\by^{\pm}| < 1\}$, choosing a branch of $\log(\by^{\pm})$ in the expression \eqref{eq:def-IG} for $\II^{\pm}_G$ specifies a single-valued section of \eqref{eq:interpret-IG}.
% \end{remark}
The proof of Theorem \ref{thm:main} is concluded by Proposition \ref{prop:UG-exists}. To understand where we are at in the proof and what still needs to be shown, the reader may wish to review the five steps of our proof as summarized at the beginning of Section \ref{sec:prooftop}; loosely speaking, we need to show that the series $\II_G^\pm$ constructed above are related by analytic continuation and a symplectic isomorphism.
\begin{proposition}\label{prop:UG-exists}
There are paths $\gamma, \delta: [0,1] \to \cM_G$ and an $F_K((z^{-1}))$-module isomorphism $\UU: H^+_G \to H^-_G$ satisfying
\begin{itemize}
\item $\gamma(0) =\delta(0) = [\epsilon:1]$ and $\gamma(1)=\delta(1)=[1:\epsilon]$ for some $\epsilon$ satisfying $0 < |\epsilon| < 1$,
\item $\UU$ is symplectic and preserves degree, and
\item $\UU\II^+_G|_{\bq = (-1)^{k-1}}$ analytically continues to $\II^-_G|_{\bq = (-1)^{k-1}}$ along $\gamma$, whereas $\UU\II^+_G|_{\bq = 1}$ analytically continues to $\II^-_G|_{\bq = 1}$ along $\delta$.
\end{itemize}
\end{proposition}

\begin{proof}
 We begin by defining $\UU$. From here, the discussion for the paths $\gamma$ and $\delta$ are in many ways parallel, but the situation giving rise to $\gamma$ is simpler.\\

\noindent
\textit{Definition of $\UU$.}
By Proposition \ref{prop:U-exists}, the isomorphism $\UU_T$ is $W$-equivariant, so in particular $\UU_T$ restricts to an isomorphism of $(H^+_T)^a$ and $(H^-_T)^a$. We write $\UU^a_T$ for this restriction, and then define $\UU$ to be the unique isomorphism of $F_K((z^{-1}))$-modules making this diagram commute:
\begin{equation}\label{eq:def-U}
\begin{tikzcd}
F_K((z^{-1})) \otimes_{R_K} A^*(X\sslash_{\theta_+}T)^a \arrow[d, "p^a", "\sim"'] \arrow[r, "\UU^a_T", "\sim"'] & F_K((z^{-1})) \otimes_{R_K} A^*(X\sslash_{\theta_-}T)^a \arrow[d, "p^a", "\sim"'] \\
F_K((z^{-1}))\otimes_{R_K} A^*(X\sslash_{\theta_+} G) \arrow[r, "\UU"] & F_K((z^{-1})) \otimes_{R_K} A^*(X\sslash_{\theta_-} G)
\end{tikzcd}
\end{equation}
Since both arrows labeled $p^a$ in \eqref{eq:def-U} shift degree by the degree of $\Delta$, we see that $\UU$ preserves degree. Similarly, since both these arrows preserve the Poincar\'e pairing up to the constant $(-1)^{|\Phi^+|}/|W|$ (Lemma \ref{lem:pairings}), we see that $\UU$ is a symplectomorphism.\\

\noindent
\textit{Proof for the path $\gamma$.}
 We define $\gamma$ to be the path appearing in Proposition \ref{prop:U-exists}. We must show that $\UU \II^+_G |_{\bq=(-1)^{k-1}}$ analytically continues to $\II^-_G|_{\bq=(-1)^{k-1}}$ along $\gamma$. We use the notation of Remark \ref{rmk:meaning}. By Lemma \ref{lem:big-I} we have
\begin{equation}\label{eq:final1}
\UU_T \II^+_T |_{q_i=1} \simeq_\gamma \II^-_T|_{q_i=1}.
\end{equation}
Observe that setting $q_i=1$ in $\II^{\pm}_T$ is the same as first setting $q_i = (-1)^{k-1}\bq$ and then setting $\bq = (-1)^{k-1}$. Moreover we may apply $j^*\partial_\Delta$ to \eqref{eq:final1}, as well to each of the $\widetilde{\II^+_{T, \ell}}$ arising in the definition of $\simeq_\gamma$, to obtain an equality 
\begin{equation}\label{eq:final2}
\left(j^* \partial_\Delta \UU_T \II^+_T |_{q_i=(-1)^{k-1}\bq} \right)|_{\bq = (-1)^{k-1}} \simeq_\gamma \left(j^*\partial_\Delta \II^-_T|_{q_i=(-1)^{k-1}\bq} \right)|_{\bq = (-1)^{k-1}}.
\end{equation}
We next commute $\partial_\Delta$ and $\UU_T$ using linearity of differentiation. Note that the symbol $\partial_\Delta$ is compatible with the implicit change of coordinates $\log(y_i^-) = - \log(y_i^+)$ by \eqref{eq:par-delta}. The specializations $y^{\pm}_i \mapsto \by^{\pm}$ are compatible with the changes of variable $y_i^+=(y_i^-)^{-1}$ and $\by^+ = (\by^-)^{-1}$, so we may further specialize both sides of \eqref{eq:final2} (as well as each $j^*\partial_\Delta \widetilde{\II^+_{T, \ell}}$ arising in the definition of $\simeq_\gamma$) to get
\begin{equation}\label{eq:final3}
\left(j^*\UU_T\partial_\Delta \II^+_T |_{q_i=(-1)^{k-1}\bq,\; y_i^+=\by^+,\; x_i = \bx_{J(i)}} \right)|_{\bq = (-1)^{k-1}} \simeq_\gamma \left(j^*\partial_\Delta \II^-_T|_{q_i=(-1)^{k-1}\bq, y^-_i = \by^-,\; x_i = \bx_{J(i)}} \right)|_{\bq = (-1)^{k-1}}.
\end{equation}
%\Wendelin{If I understand correctly, we are using here that $\partial_\Delta \II^+_T |_{q_i=(-1)^{k-1}}$ is anti-invariant and that $j^*$ induces an iso from $A^*(X//T)$ to $A^*(X//_G T)$, and use this to define a new isomorphism $\UU_T$ from $A^*(X//_G T)$ to $A^*(X//_G T)$ which we also denote $\UU_T$? In other words, we are adding a middle line in diagram (34)?} \Rachel{we could do that but let's avoid adding another flavor of $\UU$. I tried to revise.}\Wendelin{I agree, it looks good now!} 
Finally we claim that each $\partial_\Delta \widetilde{\II^+_{T, \ell}}|_{y_i^+=\by^+,\; x_i = \bx_{J(i)}}$ arising in the definition of $\simeq_\gamma$ is $W$ anti-invariant. This is true for $i=0$ by Theorem \ref{thm:nonabelianI} (since $\widetilde{\II^+_{T, 0}} = \II^+_T$), hence it is true by induction and uniqueness of analytic continuations for each $\ell$. Hence we may replace $\UU_T$ with $\UU^a_T$ in \eqref{eq:final3} and apply $(g^*)^{-1} \circ (\cup \Delta)^{-1}$ to get that
\[
\UU \II^+_G |_{\bq = (-1)^{k-1}} \simeq_\gamma \II^-_G|_{\bq = (-1)^{k-1}}.
\]

\noindent
\textit{Proof for the path $\delta$.}

We mimic the above discussion, but incorporate a change of variables. Specifically we apply $\log(y_c^-) = \log(y_c^-) - i\pi(k-1)$ to \eqref{eq:final1} followed by $j^*\partial_\Delta$. Note that the change of variables is equivalent to setting $\log(y_c^+) = \log(y_c^+) + i\pi(k-1)$. Hence if $\delta'$ is the image of $\gamma$ under this change of variables, we obtain
\[
j^*\UU_T \partial_{\Delta} \II^+_T |_{q_i=1, \;\log(y_c^+) +=i\pi(k-1)} \simeq_{\delta'} j^*\partial_\Delta \II^-_T|_{q_i=1,\; \log(y_c^-) -= i\pi(k-1)},
\]
where we have abbreviated the change of variables $\log(y_c^-) = \log(y_c^-) - i\pi(k-1)$ by $\log(y_c^-) -=  i\pi(k-1)$ and analogously for $\log(y_c^+)$.
Recall that each coordinate of $\gamma$ is the image of the path $\gamma_{\log} = \gamma_{\log}(t)$ in the $\log(y^+)$-plane pictured in \cite[Fig 1]{coates2018crepant}. This path passes between $i\pi(n-2)$ and $i\pi n$, and for $t \in [0, 1]$ near 0 or 1 satisfies $\Im(\gamma_{\log}(t))=0$ (where $\Im(z)$ denotes the imaginary part of $z$). Hence each coordinate of $\delta'$ is the image of the path $\delta'_{\log}$ in the $\log(y^+)$-plane that passes between $i\pi(n-k-1)$ and $i\pi(n-k+1)$ and for $t$ near 0 or 1 satisfies $\Im(\delta_{\log}'(t))=-\pi(k-1)$.

Consider the paths $\delta_0, \delta_1: [0,1] \to \CC_{\log(y^+)}$ given by
\[
\delta_0(t) = \log|\epsilon| - i \pi (k-1)t \quad \quad \quad \delta_1(t) = -\log|\epsilon| - i \pi (k-1)(1-t)
\]
and define
\[\delta_{\log} = \delta_1 \star \delta_{\log}' \star \delta_0.\]
 This path still passes between $i\pi(n-k-1)$ and $i\pi(n-k+1)$ but for $t$ near 0 or 1 we have that $\Im(\delta_{\log}(t))$ is close to 0.

We define $\delta$ to be the image of $(\delta_{\log}, \ldots, \delta_{\log})$ in $\cM_G$. A direct computation shows that $e^{i \pi (k-1)\sum_{j=1}^k H_j} \II^+_T$ analytically continues to $\II^+_T$ along the image of $(\delta_0, \ldots, \delta_0)$, and similarly that $\II^-_T$ analytically continues to $e^{-i \pi (k-1)\sum_{j=1}^k H_j} \II^-_T$ along the image of $(\delta_1, \ldots, \delta_1)$ (remembering to rewrite the path in terms of $y_-$). It follows that
\[
j^*\UU_T \partial_{\Delta}e^{i\pi (k-1)\sum_{j=1}^k H_j} \II^+_T |_{q_i=1, \;\log(y_c^+) += i\pi(k-1)} \simeq_{\delta} j^*\partial_\Delta e^{-i\pi(k-1)\sum_{j=1}^k H_j}\II^-_T|_{q_i=1,\; \log(y_c^-) -= i\pi(k-1)}.
\]
 We conclude using Theorem \ref{thm:nonabelianI}, specifically the formula $\eqref{eq:def IG2}$ combined with the observation \eqref{eq:simplify-sign}, and arguing as when we considered the path $\gamma$.

\end{proof}

Lemma \ref{lem:pairings} was used in the proof of Proposition \ref{prop:UG-exists}. We note that the analogous statement in cohomology when $X\sslash_\theta G$ and $X\sslash_\theta T$ are compact follows from Martin's integration formula \cite[Thm B]{martin}.

\begin{lemma}\label{lem:pairings}
%If the $K$-action on $X\sslash T$ has isolated fixed points and if the natural morphism of fixed loci $(X\sslash_G T)^K \to (X\sslash G)^K$ is surjective, then 
For any classes $\alpha, \beta \in A^*_K(X\sslash_{\pm\theta} G)_\QQ$, we have
\[
(\alpha, \beta)_{X\sslash_{\theta_{\pm}}G} = \frac{(-1)^{|\Phi^+|}}{|W|} ((p^a)^{-1}(\alpha), (p^a)^{-1}(\beta))_{X\sslash_{\theta_\pm}T}.
\] 
\end{lemma}
\begin{proof}
We begin by noting three properties of the Grassmannian flop that we use to prove the lemma. These can all be verified by direct computation. First, the $K$-actions on $X\sslash_{\pm} T$ and $X\sslash_{\pm} G$ have isolated fixed points, and the natural morphism of fixed loci $(X\sslash_{\theta_{\pm}, G} T)^K \to (X\sslash_{\theta_{\pm}} G)^K$ is surjective and all fibers have cardinality $|W|$.
% Second, from the commuting diagram of isomorphisms
% \[
% \begin{tikzcd}
% A^*_K((X\sslash_G T)^K)^a_\QQ \arrow[d, "\sim", "\iota_*"'] & A^*_K((X\sslash T)^K)^a_\QQ \arrow[l, "\sim"] \arrow[d, "\sim", "\iota_*"'] \\
% A^*_K(X\sslash_G T)^a_\QQ & A^*_K(X\sslash T)^a_\QQ \arrow[l, "\sim"', "j^*"]
% \end{tikzcd}
% \]
% it follows that if $\gamma$ is any anti-invariant class in $A^*_K(X\sslash T)$, then 
% \[
% \beta = \iota_*\sum_{F \in (X\sslash T)^K} \frac{\iota^*_F(\beta)}{e_K(N_F)} = (j^*)^{-1} \iota_*\sum_{F \in (X\sslash_G T)^K} \frac{\iota^*_Fj^*(\beta)}{e_K(N_F)}.
% \]
% It follows that
% \begin{equation}\label{eq:point2}
% \pi_* \sum_{F \in (X\sslash T)^K \seminus (X\sslash_G T)^K}
% \end{equation}
Second, if $F \in (X\sslash_{\theta_{\pm}, G} T)^K$ is a fixed point, then since $g$ is a $G/T$-fiber bundle, we have
\[
e_K(N_F) = e_K(T_{G/T})e_K(N_{g(F)}) = (-1)^{|\Phi^+|}\iota^*_F(\Delta^2) e_K(N_{g(F)}).
\]
Third, the restriction of $\Delta$ to any fixed point $F \in (X\sslash_{\theta_{\pm}} T)^K$ that does not lie in $(X\sslash_{ \theta_\pm, G} T)^K$ is zero.

Using the first two facts, we compute
\[
(\alpha, \beta)_{X\sslash_{\theta_{\pm}} G} = \pi_*\sum_{F \in (X\sslash G)^K} \frac{\iota^*_F(\alpha \cup \beta)}{e_K(N_F)} = \frac{(-1)^{|\Phi^+|}}{|W|}\pi_*\sum_{F \in (X\sslash_G T)^K} \frac{\iota^*_F(g^*\alpha \cup g^*\beta \cup \Delta^2)}{e_K(N_F)}.
\]
The numerator of the $F$-summand is equal to
\[
\iota^*_F((g^*\alpha \cup \Delta) \cup (g^*\beta \cup \Delta)).
\]
On the one hand, this is equal to $\iota^*_F((p^a)^{-1}(\alpha) \cup (p^a)^{-1}(\beta))$. On the other hand, since $\iota^*_F$ is a ring isomorphism and $\iota^*_F \Delta=0$ when $F$ is an element of $(X\sslash_{\theta_{\pm}} T)^K$ but not $(X\sslash_{ \theta_{\pm}, G} T)^K$, it is equivalent to sum over the points of $(X\sslash_{\theta_\pm} T)^K$. That is,
\[
(\alpha, \beta)_{X\sslash_{\theta_\pm }G} = \frac{(-1)^{|\Phi^+|}}{|W|}\pi_*\sum_{F \in (X\sslash_{\theta_{\pm}} T)^K} \frac{\iota^*_F((g^*\alpha \cup \Delta) \cup (g^*\beta \cup \Delta)}{e_K(N_F)}.
\]
This proves the lemma.

\end{proof}

\appendix

\section{Big points of the Lagrangian cone}\label{sec:bigappendix}

\subsection{Givental spaces and axiomatic genus-zero theories}
We recall the notion of an axiomatic genus-zero theory defined in \cite{givental2004symplectic} (but see also \cite[Sec 3]{yplee}), recasting these ideas in the language of formal schemes used in \cite[Appendix B]{coates2009computing}.

\subsubsection{Givental space}\label{sec:giventalapp}

Let $\Lambda$ be an algebra over a ground field $k$, and let $\Alg{\Lambda}$ be the category of $\Lambda$-algebras equipped with a linear, complete, and Hausdorff topology. The morphisms in $\Alg{\Lambda}$ are continuous algebra homomorphisms. For any $R \in \Alg{\Lambda}$ we have the ideal of topologically nilpotent elements
\[R^{nilp}:= \{r \in R \mid \lim_{n \to \infty} r^n = 0\}\]
as well as the ring of convergent Laurent series
\[
R\{z\} := \{\sum_{n \geq 0} r_nz^n \mid r_n \to 0 \;\text{as}\;n \to \infty\}.
\]
We let $\sO_\Lambda$ denote the tautological presheaf of rings on $\Alg{\Lambda}^{op}$; that is, $\sO_{\Lambda}(R)=R$.

Let $H$ be a $k$-vector space with a nondegenerate symmetric bilinear form $(\cdot, \cdot)$. Fix a basis $\{\phi_\alpha\}_{\alpha=1}^N$  for $H$ with a distinguished element $\phi_1=\one$ and let $\{\phi^\alpha\}_{\alpha=1}^N$ be the dual basis. We will use the Einstein summation convention for Greek indices (e.g. in the definition of $\sH(R)$ below), summing repeated Greek indices over the range $\alpha=1, \ldots, N$.

\begin{definition}\label{def:givental-space}
The \emph{Givental space} associated to $(H, \{\phi_\alpha\}, (\cdot, \cdot))$ is the presheaf $\sH$ on $\Alg{\Lambda}^{op}$ of affine $\sO_\Lambda$-modules centered at $-z$ given by
%$\sH$ from $(\Alg{\Lambda})^{op}$ infinite-dimensional modules given by
\[
\sH(R) = \{ -z + \sum_{n \in \ZZ} r^\alpha_n \phi_\alpha z^n \mid r^\alpha_n  \in R^{nilp},\; r^\alpha_n \to 0 \;\text{as}\; |n| \to  \infty \}
\]
together with the natural restriction maps.
\end{definition}
We will often use alternative coordinates $\{t^\alpha_n, p_{\beta \ell}\}$ on $\sH(R)$ given by the identification 
%\Rachel{changed index set on $p_{\beta \ell}$}
\[
-z + \sum_{n \in \ZZ} r^\alpha_n \phi_\alpha z^n =: -z + \sum_{n \geq 0} t^\alpha_n \phi_\alpha z^n + \sum_{\ell \geq 0} p_{\beta \ell} \frac{\phi^\beta}{(-z)^{\ell + 1}}.
\]
We will abbreviate the above equality by writing $-z + \br = -z + \bt + \bp$. These coordinates identify a sub-presheaf $\sH_+ \subset \sH$ given by $\sH_+(R) = \{I \in \sH(R) \mid I = -z + \bt \}$.

One defines the tangent sheaf $T\sH$ on $\Alg{\Lambda}^{op}$ by setting $T\sH(R) := \sH(R[\epsilon]/(\epsilon^2))$ as in \cite[Appendix B]{coates2009computing}. The tangent space at an element $I$ of $\sH(R)$ is the preimage of $I$ under the map $T\sH \to \sH$ induced by $R[\epsilon]/(\epsilon^2) \to R$; explicitly, is given by
\begin{align}
T_I\sH &= \{-z + I + \epsilon\sum_{n \in \ZZ} \dot{r}^\alpha_n \phi_\alpha z^n \mid \dot{r}^\alpha_n \to 0 \;\text{as}\; |n| \to \infty \} \notag\\
&\simeq  \{\sum_{n \in \ZZ} \dot{r}^\alpha_n \phi_\alpha z^n \mid \dot{r}^\alpha_n \to 0 \;\text{as}\; |n| \to \infty \}. \label{eq:tangent-iso}
\end{align}
To clarify, a priori the fiber $T_I\sH$ is an affine $R$-module centered at $-z+I$, but via the above isomorphism we view it as an $R$-module and even as an $R\{z\}$-module. 
We will often use alternative coordinates on $T_I\sH$, writing 
%\Rachel{changed index set on $\dot{p}$}\Wendelin{I don't remember what it was before, but it looks good to me. }
\[
\sum_{n \in \ZZ} \dot{r}^\alpha_n \phi_\alpha z^n =: \sum_{n \geq 0} \dot{t}^\alpha_n \phi_\alpha z^n + \sum_{\ell \geq 0} \dot{p}_{\beta \ell} \frac{\phi^\beta}{(-z)^{\ell+1}}
\]
or $\bdr = \bdt + \bdp$ for short.

\subsubsection{Axiomatic genus-zero theories}
\label{sec:axiomatic g0}

%Let $\bullet$ denote the functor from $\Alg{\Lambda}$ to sets that sends a $\Lambda$-algebra $R$ to itself (as a set). 
A \emph{function on $\sH_+$} is a morphism of presheaves on $\Alg{\Lambda}^{op}$  from $\sH_+$ to $\sO_\Lambda$, denoted by
\[
-z + \bt \mapsto \llangle \rrangle_{\bt}.
\]
If $m$ is a positive integer and $R \in \Alg{\Lambda}^{op}$, then $ R\llbracket s_1, \ldots, s_m \rrbracket$ is also in $\Alg{\Lambda}^{op}$.
% Observe that if $m$ is a positive integer and $R \in \Alg{\Lambda}^{op}$, then $R\llbracket \bs \rrbracket := R\llbracket s_1, \ldots, s_m \rrbracket$ is also in $\Alg{\Lambda}^{op}$, where if the topology on $R$ is induced by the system of ideals $\{I_n\}_{n \geq 0}$ then the topology on $R\llbracket \bs\rrbracket$ is induced by the system
% \[\big \{I_n\llbracket s_1, \ldots, s_m \rrbracket + (s_1, \ldots, s_m)^n R\llbracket s_1, \ldots, s_m \rrbracket \big \}_{n\geq 0}.\]
% Moreover we have well-defined endomorphisms $\partial/\partial s_i$ of $ R\llbracket \bs \rrbracket$. Using these observations, 
We define the \textit{derivative}
\[
\llangle \phi_{\beta_1} \psi^{\ell_1}, \ldots,  \phi_{\beta_m} \psi^{\ell_m}\rrangle_\bt := \left( \frac{\partial}{\partial s_1}\ldots  \frac{\partial}{\partial s_m} \llangle \rrangle_{\bt + s_1\phi_{\beta_1} \psi^{\ell_1} + \ldots + s_m \phi_{\beta_m}\psi^{\ell m}} \right) \Big |_{s_1=\ldots=s_m=0}
\]
noting that if $-z + \bt \in \sH_+(R)$ then the right hand side is naturally an element of $R$; in particular, this derivative is again a function on $\sH_+$.
 By linearity we define derivatives with respect to arbitrary elements of $H$ (not just basis elements).

 Let $\tau^1, \ldots, \tau^N$ be formal variables (corresponding to the basis $\phi_1, \ldots, \phi_N$ of $H$) and define
 \[
 \widetilde{\Lambda} := \Lambda\llbracket \tau^1, \ldots, \tau^N \rrbracket.
 \]
\begin{definition}
A \emph{genus-zero theory} is a function $\llangle \rrangle_\bt$ on $\sH_+$ that satisfies %\Rachel{changed next equation}
\begin{align}\label{eq:shape}
\llangle \rrangle_{\bt+\bs} = \llangle \rrangle_{\bt} + \sum_{k \geq 0} s^\alpha_k \llangle \phi_\alpha \psi^k \rrangle_\bt  &+ \frac{1}{2}\sum_{k_1, k_2 \geq 0}  s^{\alpha_1}_{k_1}s^{\alpha_2}_{k_2} \llangle \phi_{\alpha_1} \psi^{k_1} , \phi_{\alpha_2} \psi^{k_2}\rrangle_\bt\\
&+ \frac{1}{3!} \sum_{k_1, k_2, k_3 \geq 0}  s^{\alpha_1}_{k_1}s^{\alpha_2}_{k_2}s^{\alpha_3}_{k_3} \llangle \phi_{\alpha_1} \psi^{k_1} , \phi_{\alpha_2} \psi^{k_2}, \phi_{\alpha_3} \psi^{k_3}\rrangle_\bt + \;\ldots\notag
\end{align}
% \begin{align}\label{eq:shape}
% \llangle \rrangle_\bt = K + \sum_{k \geq 0} t^\alpha_k \llangle \phi_\alpha \psi^k \rrangle_\bt  &+ \frac{1}{2}\sum_{k_1, k_2 \geq 0}  t^{\alpha_1}_{k_1}t^{\alpha_2}_{k_2} \llangle \phi_{\alpha_1} \psi^{k_1} , \phi_{\alpha_2} \psi^{k_2}\rrangle_\bt\\
% &+ \frac{1}{3!} \sum_{k_1, k_2, k_3 \geq 0}  t^{\alpha_1}_{k_1}t^{\alpha_2}_{k_2}t^{\alpha_3}_{k_3} \llangle \phi_{\alpha_1} \psi^{k_1} , \phi_{\alpha_2} \psi^{k_2}, \phi_{\alpha_3} \psi^{k_3}\rrangle_\bt + \;\ldots\notag
% \end{align}
where the infinite sum is defined to be the limit of the partial sums in the topological ring. We moreover require that $\llangle \rrangle_\bt$ satisfy the following differential equations:
 \begin{description}
\item[(DE)] 
$\llangle \psi \rrangle_\bt = \sum_{k \geq 0} t^\alpha_k \llangle \phi_\alpha \psi^k \rrangle_\bt - 2 \llangle \rrangle_\bt$
\item[(SE)]
$\llangle \one \rrangle_\bt = \frac{1}{2}(\bt_0, \bt_0) + \sum_{k \geq 0} t^\alpha_{k+1} \llangle \phi_\alpha \psi^k \rrangle_\bt$
where $\bt_0 = \sum_\alpha t^\alpha_0 \phi_\alpha$ 
\item[(TRR)]
$\llangle \phi_\alpha\psi^{k+1} ,\phi_\beta\psi^\ell ,\phi_\gamma \psi^m\rrangle_\bt=\llangle \phi_\alpha \psi^k, \phi_\nu \rrangle_\bt \llangle \phi^\nu, \phi_\beta\psi^\ell ,\phi_\gamma \psi^m\rrangle_\bt$
\end{description}
Finally, we require that the \emph{J-function} 
%\Rachel{changed index set on $\ell$}
\begin{equation}\label{eq:Jfunc}
J := -z + \tau^\alpha \phi_\alpha + \sum_{\ell \geq 0} \llangle \phi_\beta \psi^\ell \rrangle_{\bt} \frac{\phi^\beta}{(-z)^{\ell+1}}
\end{equation}
be an element of $\sH(\widetilde{\Lambda})$, and that 
\begin{equation}\label{eq:Jisbig}
\frac{\partial J}{\partial \tau^\alpha} \equiv \phi_\alpha \quad \quad \text{modulo} \quad \widetilde{\Lambda}^{nilp}.
\end{equation}
\end{definition}
\begin{remark}
 The genus-zero potential in Gromov-Witten theory defines an axiomatic genus-zero theory. Equation \eqref{eq:shape} follows from the linearity property of Gromov-Witten correlators. The assertions that the $J$-function \eqref{eq:Jfunc} is an element of the Givental space and that \eqref{eq:Jisbig} holds follow from the vanishing of degree-zero correlators with 1 or 2 insertions. 
\end{remark}

We define the \emph{Givental cone} associated to a genus-zero theory to be the subsheaf $\sL \subset \sH$ on $\Alg{\Lambda}^{op}$ given by 
\[
\sL(R) = \{-z + \bt + \bp \in \sH(R) \mid p_{k\alpha} = \llangle \phi_\alpha \psi^k \rrangle_\bt\}.
\]
% One defines the tangent functor $T\sL$ by $T\sL(R) = \sL(R[\epsilon]/(\epsilon^2)$ and for $I \in \sL$ the tangent space to be the fiber of $T\sL \to \sL$. One checks that with these definitions, if $I = -z + \bt + \bp$, we have \Rachel{is this true?? I didn't actually check!}
We define $T\sL$ and $T_I\sL$ in analogy with $T\sH$ and $T_I\sH$. One checks using \eqref{eq:shape} that under the isomorphism \eqref{eq:tangent-iso}, we have 
%\Qaasim{Is $\llangle\rrangle_{\bt + \epsilon \bdt} = \llangle\rrangle_{\bt} + \epsilon \llangle\rrangle_{\bdt}$?}\Rachel{not quite. details are in a hidden comment below the def of $T_I\sL$. I don't really want to write them in here because I think my notation might be off and I don't want to fix it, but feel free to uncomment and take a look.}
\begin{equation}\label{eq:cone-tangents}
T_I\sL = \left\{ \bdt + \bdp \in T_I\sH \mid \dot{p}_{\ell\beta}=\sum_{\alpha, k \geq 0}\dot{t}_{k}^\alpha \llangle \phi_\alpha \psi^k, \phi_\beta \psi^\ell \rrangle_{\bt}\right\}.
\end{equation}
%%%% DETAILS ON THE ABOVE:
% If we write $-z + \bt + \epsilon \dot{\bt} + \bp + \epsilon \dot{\bp} $ then this is a point of $\sL$ if and only if
% \[
% p_{\ell \beta} + \epsilon \dot{p}_{\ell \beta} = \llangle \phi_\beta \psi^\ell \rrangle_{\bt + \epsilon \dot{\bt}} = \llangle \phi_\beta \psi^\ell \rrangle_{\bt } + \epsilon \sum \dot{t}^\alpha_k \llangle \phi_\alpha \psi^k, \phi_\beta \psi^\ell \rrangle_{\bt}
% \]
 %using \eqref{eq:shape} and the fact that $\epsilon^2=0$. Now equate $\epsilon$ parts and non-$\epsilon$ parts.

\begin{remark}
The results in \cite[Appendix B]{coates2009computing} hold for all axiomatic genus-zero theories. 
%For example . . . \Rachel{We should  check the citation of Dijkgraaf-Witten on p.44. I did not look it up.}
\end{remark}

\subsection{Big points of the Givental cone}
Let $(H, (\cdot, \cdot), \{\phi_\alpha\})$ be as in Section \ref{sec:giventalapp} and let $\sH$ be the associated genus-zero theory. 

\begin{definition}\label{def:big} An element $I$ of $\sH(\widetilde{\Lambda})$ is \emph{big} if
\begin{equation}\label{eq:big-deriv}
\frac{\partial I}{\partial \tau^\alpha} \equiv \phi_\alpha \quad\quad\text{modulo}\quad \widetilde{\Lambda}^{nilp}
\end{equation}
\end{definition}

The significance of big points is that they completely determine the Givental cone. To state this precisely (Proposition \ref{prop:bigI} below), we introduce some notation. For any $I = -z + \bt_I + \bp_I$ in $\sH(R)$, we define 
\[
\tau_{I}^{\alpha} := \llangle \one, \phi^\alpha \rrangle_{\bt_I} \in R.
\]
The elements $\tau_{I}^{\alpha} \in R$ induce a topological $\Lambda$-algebra homomorphism $\tau_I: \widetilde{\Lambda} \to R$. For any big point $I$, we define an $N \times N$ matrix $DI$ with entries $(DI)_{\alpha \beta}$ in $\widetilde{\Lambda}\{z, z^{-1}\}$ 
%\Rachel{is this the correct place for $DI$ to live? also multiplication by $w$ is well defined because you can look modulo any ideal defining the topology of $R$, right? } 
given by
\[
(DI)_{\alpha \beta} = \frac{\partial I^{\alpha}}{\partial \tau^\beta} \quad \quad \quad \text{where} \;\; I = \sum_{\alpha=1}^N I^\alpha \phi_\alpha.
\]
Finally, if $\sigma: S \to R$ is a homomorphism of topological $\Lambda$-algebras and $f$ is an element of $S\{z\}$, we denote by $f(\sigma)$ the element of $R\{z\}$ obtained by applying $\sigma$ to the coefficients of $f$.
\begin{proposition}\label{prop:bigI}
Let $I$ be a big point of the Givental cone. Then $\tau_I$ is an automorphism of $\widetilde{\Lambda}$, and an element $K$ of $\sH(R)$ is in the Givental cone if and only if
\begin{equation}\label{eq:a}
K = zDI(t) w
\end{equation}
for some topological $\Lambda$-algebra homomorphism $t: \widetilde{\Lambda} \to R$ and $w \in R\{z\}^{N}$ satisfying $w=(-1, 0, \ldots, 0)$ modulo $R^{nilp}$. In this case $t$ and $w$ are unique and $t = \tau_K \circ \tau_I^{-1}$.
\end{proposition}
\begin{remark}\label{rmk:format}
Any topological $\Lambda$-algebra homomorphism $t: \widetilde{\Lambda} \to R$ corresponds to elements $t^1, \ldots, t^N \in R$. The fact that our homomorphisms are continuous implies that the $t^\alpha$ are in fact elements of $R^{nilp}$. 

Moreover, any topological $\Lambda$-algebra homomorphism $t: \widetilde{\Lambda} \to R$ and $w \in R\{z\}^N$ satisfying \eqref{eq:a} will also satisfy $w \equiv -\one $ modulo $R^{nilp}$. This follows from reducing \eqref{eq:a} modulo $R^{nilp}$: the left hand side becomes $-z$ while the right hand side becomes $zw$ since $I$ is big.

Conversely, if $t: \widetilde{\Lambda} \to R$ is a topological $\Lambda$-algebra homomorphism and $w \in R\{z\}^N$ satisfies $w \equiv -\one $ modulo $R^{nilp}$, then $zDI(t)w$ will be an element of $\sH(R)$. 
\end{remark}
\begin{example}
The $J$-function \eqref{eq:Jfunc} is a big point of the Givental cone of $\sH$, and the homomorphism $\tau_J$ is just the identity. Hence Proposition \ref{prop:bigI} says that $K \in \sH(R)$ is on the Givental cone if and only if $K=zDJ(\tau_K) w$ for some $w \in R\{z\}^N$. See also Lemma \ref{lem:Jcase} below.
\end{example}

The proof of Proposition \ref{prop:bigI} was explained to us by Hiroshi Iritani. We prove it in the special case when $I=J$ before proceeding with the proof of the general statement. Let $[\cdot]_+$ denote truncation of the argument to nonnegative powers of $z$.

\begin{lemma}\label{lem:uniqueness}
The assignment $(t, w) \mapsto [z+zDJ(t)w]_+$ 
%\Qaasim{$+z$?}\Wendelin{I agree that this should be $+z$ going by Hiroshi's email.} 
is injective, where $t: \widetilde{\Lambda} \to R$ is a homomorphism and $w \in R\{z\}^N$ satisfies $w \equiv -\one$ modulo $R^{nilp}$.
%\Qaasim{$R^{nilp}$?}.
\end{lemma}
\begin{proof}
Let $\bt = [z + zDJ(t)w]_+$ and write $z + zDJ(t)w = \bt + \bp$ for some $\bp$. We have
\begin{equation}\label{eq:u1}
zw = DJ(t)^{-1} (-z + \bt + \bp) = [DJ(t)^{-1}(-z + \bt)]_+.
\end{equation}
%\Qaasim{I don't currently understand the second equality}\Wendelin{I think we are taking $+$ on both sides of the first equality, and using that $zw$ only has positive terms in z and $DJ(t)^{-1}\bp$ only has negative terms in $z$}\Qaasim{Thanks!}
From here it is clear that if we can show $t$ is determined by $\bt$, then $w$ is also determined. Taking the $z$-constant terms of \eqref{eq:u1} gives
\[
0 = [DJ(t)^{-1}(-z + \bt)]_0
\]
where $[\cdot ]_0$ denotes the constant term. If we expand $DJ(t)^{-1} = Id + D_1(t) z^{-1} + D_2(t) z^{-2} + \ldots$ and $\bt = \bt_0 + \bt_1z + \bt_2 z^2 + \ldots$, then this equation is equivalent to
\begin{equation}\label{eq:u2}
0 = -D_1(t)\one + \bt_0 + D_1(t)\bt_1 + D_2(t)\bt_2 + \ldots.
\end{equation}
We claim $-D_1(t)\one=t$; indeed, we have $DJ(t) = Id - D_1(t)z^{-1} + O(z^{-2})$, so $-D_1(t)\one$ is the coefficient of $z^{-1}$ in $\partial J/\partial \tau^1(t)$. This coefficient is $t$ by the string equation.

Continuing on, let $t$ and $t'$ be two solutions to \eqref{eq:u2} and let $I \subset R$ be an open ideal. It is enough to show that $t-t'$ is in $I$. Since $\bt_k \in I$ for sufficiently large $k$, there is some $k$ such that
\begin{equation}\label{eq:inductme}
t - t' \equiv -(D_1(t) - D_1(t'))\bt_1 - (D_2(t) - D_2(t'))\bt_2 - \ldots - (D_k(t) - D_k(t'))\bt_k\quad \quad \text{modulo} \;\; I.
\end{equation}
Here we computed $t-t'$ using \eqref{eq:u2} and the equalities $-D_1(t) \one = t$ and $-D_1(t') \one  = t'.$
Now let $J \subset R$ be the ideal generated by the coordinates of $t-t'$ and $\bt_1, \bt_2, \ldots, \bt_k$; in particular, $J$ is generated by finitely many nilpotent elements (see Remark \ref{rmk:format}). Hence there is an integer $n$ for which $J^n \subset I$, and if we can show $t-t'$ is in $J^n$ we are finished.

We show $t-t' \in J^n $ for all $n$ by induction. The base case $n=1$ holds by construction of $J$. If $t-t' \in J^n$, then $D_i(t) - D_i(t') \in J^n$ for all $i$: this is because we may write the telescoping sum 
%\Rachel{removed a minus sign from the right hand side. can other people verify that it is correct now??}\Wendelin{Looks good to me}
\begin{align*}
D_i(t) - D_i(t') = \sum_{\ell=1}^N D_i(t_1, \ldots, t_{\ell-1}, t_\ell, t_{\ell+1}', \ldots, t_N') - D_i(t_1, \ldots, t_{\ell-1}, t_\ell', t_{\ell+1}', \ldots, t_N')
\end{align*}
and each entry of the $\ell^{th}$ summand may be written as
\[
p(t_\ell) - p(t_\ell') = \sum_{j=1}^M a_j\cdot(t_\ell^j - t_\ell'^j) = (t_\ell-t_\ell')\sum_{j=1}^M a_j\cdot (t_\ell^{j-1} + \ldots + t_\ell'^{j-1})
\]
for some univariate polynomial $p$ of degree $M$ with coefficients $a_j$ in $\Lambda[t_1, \ldots, t_{\ell-1}, t_{\ell+1}', \ldots, t_N']$. It follows that the right hand side of \eqref{eq:inductme} is in $J^{n+1}$, hence the left hand side is too.
%\Wendelin{I think all of the $J^n+I$ in the above proof should just be $J^n$.}
\end{proof}

\begin{lemma}\label{lem:Jcase}
Proposition \ref{prop:bigI} holds for $I=J$.
\end{lemma}
\begin{proof}
Let $K \in \sH(R)$ be a point of the Givental cone. We show that $K = zDJ(\tau_K)w$ for some $w\in R\{z\}^N$. It is enough to show that $z^{-1}K$ is an element of $T_K\sL$; then \cite[Prop B.4]{coates2009computing} gives the required $w$. For this, let $K=-z + \bt + \bp$. Then $z^{-1}K = -1 + z^{-1}(\bt + \bp)$ may be viewed as the element $\bdt + \bdp$ of $T_K\sH$ whose coordinates are given by
\begin{align}
\dot{t}_\ell^\alpha&:=t_{\ell+1}^\alpha &(\text{for}\; (\ell, \alpha) \neq (0, 1) \label{eq:shift-coordinates}\\
\dot{t}_0^1&:=-1+t_{1}^1 \notag\\
\dot{p}_{0\alpha}&:=\sum_\beta-t_{0}^\beta g_{\beta \alpha}\notag\\
\dot{p}_{\ell\alpha}&:=-p_{({\ell-1})\alpha}  & (\text{for}\; \ell > 0).\notag
\end{align}
where $g_{\beta \alpha}:= (\phi_\beta, \phi_\alpha)$.
To show that $z^{-1}K \in T_K\cL$ it is enough to show that the equation
\begin{equation}\label{eq:stringproblem}
\dot{p}_{\ell\beta}=\sum_{\alpha, k \geq 0}\dot{t}_{k}^\alpha \llangle \phi_\alpha \psi^k, \phi_\beta \psi^\ell \rrangle_{\bt}.
\end{equation}
is satisfied (see \eqref{eq:cone-tangents}).
In the case $\ell > 0$, we use $p_{\ell\beta}=\llangle \phi_\beta \psi^l \rrangle_{\bt}$ to obtain that \eqref{eq:stringproblem} is equivalent to
\[
-\llangle \phi_\beta \psi^{\ell-1} \rrangle_{\bt}=\sum_{\alpha, k \geq 0}t_{k+1}^\alpha \llangle \phi_\alpha \psi^k, \phi_\beta \psi^\ell \rrangle_{\bt}-\llangle 1, \phi_\beta \psi^\ell \rrangle_{\bt} 
\]
Rearranging, this agrees with the derivative of the string equation (SE) in the $\phi_\beta \psi^\ell$-direction,\footnote{By this, we mean to evaluate both sides of (SE) at $\bt + s \phi_\beta \psi^\ell$, differentiate both sides with respect to $s$, and then set $s=0$.} where we use the product rule to differentiate terms $t^\alpha_{k+1}\llangle \phi_\alpha \psi^k \rrangle_\bt$ with $k+1=1$ and $\alpha=\beta$.
In the case $\ell=0$, equation \eqref{eq:stringproblem} is equivalent to 
\[
-\sum_\beta t_{0}^\beta g_{\beta \alpha}=\sum_{\alpha, k \geq 0}t_{k+1}^\alpha \llangle \phi_\alpha \psi^k, \phi_\beta\rrangle_{\bt}-\llangle 1, \phi_\beta \rrangle_{\bt}
\]
which again agrees with the derivative of the string equation in the $\phi^\beta$-direction.

Conversely, let $K \in \sH(R)$ be an element satisfying $K = zDJ(t)w$ for some pair $(t, w)$. If we write $K = -z + \bt_K + \bp_K$, then there is some element $L = -z + \bt_K + \bp \in \sH(R)$ that is on the cone. But $L = zDJ(t')w'$ for some pair $(t', w')$ by the previous direction, and
\[
[z + zDJ(t')w']_+ = \bt_K = [z+zDJ(t)w ]_+.
\]
It follows from Lemma \ref{lem:uniqueness} that $(t, w) = (t', w')$ and hence $K=L$ is on the cone.

The uniqueness of the pair $(t, w)$ in Proposition \ref{prop:bigI} also follows from Lemma \ref{lem:uniqueness}.

\end{proof}

\begin{lemma}\label{lem:technical}If $ I \in \sH(\widetilde{\Lambda})$ is a big point of $\sL$ then $\tau_I$ is an automorphism of $\widetilde{\Lambda}$ and 
\[
DI = DJ(\tau_I) V
\]
for some $N\times N$ matrix $V$ with entries in $\widetilde{\Lambda}\{z\}$ such that $V$ reduces to the identity modulo $\widetilde{\Lambda}^{nilp}$; in particular, $V$ is invertible as an endomorphism of $(\widetilde{\Lambda}\{z\})^N$.
\end{lemma}

\begin{proof}
By Lemma \ref{lem:Jcase} we may write
\[
I = zDJ(\tau_I)w =  z\sum_\alpha \frac{\partial J}{\partial \tau^\alpha}(\tau_I) w^\alpha
\]
for some $w =w^\alpha \phi_\alpha $ with $w^\alpha \widetilde{\Lambda}\{z\}$ satisfying $ w\equiv -\one$ modulo $\widetilde{\Lambda}^{nilp}$. 
Differentiation with respect to $\tau_\beta$ yields
\begin{equation}\label{eq:thing1}
\frac{\partial I}{\partial \tau^\beta} =  \sum_\alpha z\left(\sum_\gamma \frac{\partial J}{\partial \tau^\gamma \partial \tau^\alpha}(\tau_I)\frac{\partial \tau^\gamma_I}{\partial \tau^\beta} w^\alpha \right)+ z\frac{\partial J}{\partial \tau^\alpha}(\tau_I) \frac{\partial w^\alpha}{\partial \tau^\beta}.
\end{equation}
It follows from the topological recursion relation that 
%\Rachel{corrected the sign below and later. please check that you agree!}\Wendelin{agree!}
\[
-z \frac{\partial J}{\partial \tau^\gamma \partial \tau^\alpha}(\tau_I) = \llangle \phi_\gamma, \phi_\alpha, \phi^\nu \rrangle_{\tau_I^\epsilon \phi_\epsilon} \frac{\partial J}{\partial \tau^\nu}(\tau_I).
\]
This equality is also known as the quantum differential equation. Combining with \eqref{eq:thing1} yields
\begin{equation}\label{eq:thing2}
\frac{\partial I}{\partial \tau^\beta} =  \sum_\alpha \left(\sum_\gamma \sum_\nu -\llangle \phi_\gamma, \phi_\alpha, \phi^\nu \rrangle_{\tau_I^\epsilon \phi_\epsilon} \frac{\partial J}{\partial \tau^\nu}(\tau_I) \frac{\partial \tau^\gamma_I}{\partial \tau^\beta} w^\alpha \right)+ z\frac{\partial J}{\partial \tau^\alpha}(\tau_I) \frac{\partial w^\alpha}{\partial \tau^\beta} = DJ(\tau_I)V_\beta
\end{equation}
where $V_\beta$ is the $\beta^{th}$ column of the $N \times N$ matrix $V$ given by
\[
V_{\nu \beta} = \sum_\alpha \sum_\gamma -\llangle \phi_\gamma, \phi_\alpha, \phi^\nu \rrangle_{\tau_I^\epsilon \phi_\epsilon}  \frac{\partial \tau^\gamma_I}{\partial \tau^\beta} w^\alpha + z\frac{\partial w^\nu}{\partial \tau^\beta} \in \widetilde{\Lambda}\{z\}.
\]
We note that $V$ reduces to the identity modulo $\widetilde{\Lambda}^{nilp}$ because both $DI$ and $DJ(\tau_I)$ do.

On the other hand, reducing \eqref{eq:thing2} modulo $\widetilde{\Lambda}^{nilp}$ yields
\begin{equation}\label{eq:thing3}
\delta_{\epsilon \beta} \equiv \frac{\partial \tau^\epsilon_I}{\partial \tau^\beta} + z\frac{\partial w^\epsilon}{\partial \tau^\beta} \quad \quad \text{modulo}\;\;\widetilde{\Lambda}^{nilp}.
\end{equation}
To obtain the first summand on the right hand side we used the following, in this order: First, $\tau_I$ preserves $\widetilde{\Lambda}^{nilp}$, so $({\partial J^\epsilon}/{\partial \tau^\nu})(\tau_I)$ is equal to $\delta_{\epsilon \nu}$ modulo $\widetilde{\Lambda}^{nilp}$ since $J$ is a big point; %(\Wendelin{and similarly for $I$}; 
second, $w^\alpha \equiv -\delta_{1\alpha}$ modulo $\widetilde{\Lambda}^{nilp}$; third, $\llangle \phi_\gamma, \one, \phi^\nu\rrangle_{\tau^\epsilon_I \phi_\epsilon} \equiv (\phi_\gamma, \phi^\nu) = \delta_{\gamma \nu}$ modulo $\widetilde{\Lambda}^{nilp}$ by the string equation. %\Rachel{please check this entire computation carefully!} \Wendelin{Looks good!}\Qaasim{confused about where the minus sign due to $w^{\alpha}$ goes?}
Comparing coefficients of $z$ in \eqref{eq:thing3} we find
$
 {\partial \tau^\epsilon_I}/{\partial \tau^\beta}\equiv \delta_{\epsilon \beta} $ modulo $\widetilde{\Lambda}^{nilp},$
showing that the homomorphism $\tau_I: \widetilde{\Lambda} \to \widetilde{\Lambda}$ is invertible.

\end{proof}

\begin{proof}[Proof of Proposition \ref{prop:bigI}] Let $K \in \sH(R)$. By Lemma \ref{lem:Jcase} we have that $K$ is in the Givental cone if and only if 
\begin{equation}\label{eq:ug}
K = zDJ(t)v
\end{equation}
for some (unique) pair $(t, v)$. By Lemma \ref{lem:technical} we may write $DJ(t) = DI(t\circ \tau_I^{-1}) V'$ where $V' = V(t \circ \tau_I^{-1})$. So if $K$ is on the Givental cone then by Lemma \ref{lem:Jcase} we have $t = \tau_K$ in \eqref{eq:ug}, and
\[
K = z DI(\tau_K \circ \tau_I^{-1}) V'v
\]
and we obtain one direction of the lemma by setting $w = V'v$. Conversely if $K = zDI(t)w$ for some pair $(t, w)$, then replacing $DI(t)$ by $DJ(t \circ \tau_I) V(t)$ we see that $K$ is on the cone. 

Finally, for uniqueness of the pair $(t, w)$ in Proposition \ref{prop:bigI}, suppose that $K \in \sH(R)$ satisfies
\[
K = zDI(t)w = zDI(t')w'.
\]
By Lemma \ref{lem:technical} we have
\[
zDJ(t \circ \tau_I) V(t)w = zDJ(t' \circ \tau_I) V(t')w'
\]
but now $t=t'$ and $w=w'$ follows from the fact that Proposition \ref{prop:bigI} holds when $I=J$ (or more specifically Lemma \ref{lem:uniqueness}).
\end{proof}
\bibliographystyle{abbrv}
\bibliography{references}

\end{document}